\newif\ifShortTable
\ShortTabletrue

\newif\iftcnsVersion
\tcnsVersiontrue

\newif\ifOneColumn
\OneColumnfalse

\ifOneColumn
\documentclass[journal, onecolumn,usenames,dvipsnames]{IEEEtran}
	\renewcommand{\baselinestretch}{2}
\else
\documentclass[journal,usenames,dvipsnames]{IEEEtran}
\fi

\newif\ifArxivVersion
\ArxivVersiontrue

\newif\ifNotArxivVersion

\usepackage{etex}
\reserveinserts{28}

\usepackage{makeidx}         
\usepackage{graphicx}        
\usepackage{multicol}        

\makeindex             

\usepackage{enumitem, fancyhdr, amsthm, amsmath, amsfonts,  amssymb, csquotes, mathtools, url, graphicx, balance, comment, caption, mathrsfs, upgreek, wrapfig, cases, pgfplots, subfig, caption, fancybox,  multirow, hhline, mathtools, chngcntr, booktabs, algpseudocode, mathabx, floatrow, algorithm, algorithmicx, mathrsfs, etex, etoolbox, xfrac,  cite, xstring, xparse, ifmtarg, xifthen, bm, array, multicol, layouts, printlen,  soul, cancel, tabularx, ragged2e, relsize, scalefnt, varwidth, setspace}

\algnewcommand\algorithmicforeach{\textbf{for each}}
\algdef{S}[FOR]{ForEach}[1]{\algorithmicforeach\ #1\ \algorithmicdo}

\newcolumntype{x}[1]{>{\centering\let\newline\\\arraybackslash\hspace{0pt}\vspace{1pt}}p{#1}}
 
\usepackage[disable]{todonotes}
\usepackage[numbers]{natbib}
\usepackage[normalem]{ulem}
\pgfplotsset{compat=1.13}

\usepackage{tikz,lmodern}
\usetikzlibrary{patterns,arrows.meta,shapes,snakes,automata,backgrounds,petri,fit,calc,positioning,external,intersections,trees}
\usepackage{color,xcolor}
\usepackage[colorlinks=true,
            raiselinks=true,
            linkcolor=MidnightBlue,
            citecolor=Mahogany,
            urlcolor=ForestGreen,
            pdftitle={Resilience of Electricity Distribution Networks},
            pdfkeywords={DN Security},
            pdfsubject={Initial Manuscript},
            plainpages=false
            ]{hyperref}

\newcommand{\xfnm}[1][]{\ifx!#1!\else\unskip,\space#1\fi} 
\usepackage[bottom,multiple]{footmisc} 

\newenvironment{condition}[1][]{\setcounter{subcondcounter}{-1}\hspace{-0.2cm}\vspace{-0.15cm}\refstepcounter{condcounter} 
	\par\medskip
	\hspace{-0.4cm}\textbf{(C\thecondcounter #1)} 
}{
}

\newcounter{algsubstate}


\newcommand{\subparagraph}{}

\usepackage[compact]{titlesec}

\usepackage[nameinlink]{cleveref}

\newif\ifLongVersion
\LongVersionfalse

\ifArxivVersion

	\NotArxivVersionfalse
\else
	
	\renewcommand{\baselinestretch}{1}

	\NotArxivVersiontrue
\fi

\crefname{section}{Sec.}{§§}
\Crefname{section}{Sec.}{§§}
\crefname{subsection}{Sec.}{§§}
\crefname{equation}{Eq.}{}
\Crefname{equation}{Equation}{Equations}

\newtheoremstyle{theoremdd}
{5pt}
{5pt}
{\itshape}
{0pt}
{\bfseries}
{.}
{ }
{\thmname{#1}\thmnumber{ #2}\thmnote{ (#3)}}

\theoremstyle{theoremdd}

\newtheorem*{formulation*}{}

\theoremstyle{remark}
\newtheorem{remark}{Remark}

\crefname{table}{Table}{Tables}
\crefname{figure}{Fig.}{Figures}
\crefname{theorem}{Theorem}{Theorems}
\crefname{proposition}{Proposition}{Propositions}
\crefname{lemma}{Lemma}{Lemmas}
\crefname{algorithm}{Algorithm}{Algorithms}
\crefname{myclaim}{Claim}{Claims}

\makeatletter
\def\thm@space@setup{%
  \thm@preskip=0.2cm plus 0.2cm minus 0.2cm
  \thm@postskip=\thm@preskip 
}
\renewenvironment{proof}[1][\proofname]{\par
  \pushQED{\qed}%
  \normalfont
  \topsep0pt \partopsep5pt 
  \trivlist
  \item[\hskip\labelsep
        \itshape
    #1\@addpunct{.}]\ignorespaces
}{%
  \popQED\endtrivlist\@endpefalse
  \addvspace{6pt plus 6pt} 
}
\makeatother

\makeatletter

\newtheoremstyle{named}{}{}{\itshape}{}{\bfseries}{.}{.5em}{\thmnote{#3's }#1}
\theoremstyle{named}

\newtheoremstyle{mynamed}{}{}{\itshape}{}{\bfseries}{.}{.5em}{#1 \thmnote{A }}
\theoremstyle{named}

\newcommand{\abs}[1]{\left\lvert{#1}\right\rvert}

\allowdisplaybreaks

\counterwithout{equation}{section} 
\counterwithout{figure}{section} 

\algrenewcommand{\algorithmiccomment}[1]{\hskip3em$\slash\slash$ #1}

\newcommand{\linfinityNorm}[1]{\left\lVert#1\right\rVert_\infty}

\def \unity {\mathbf{1}}
\def \zero {\mathbf{0}}

\def \setDER {\mathcal{S}}

\def \setLoadControl {\mathcal{B}}
\def \setMicrogrid {\mathcal{M}}

\newcommand{\N}[1][]{\mathcal{N}_{#1}}
\newcommand{\E}[1][]{\mathcal{E}_{#1}}

\def \M {\mathcal{M}}
\def \bigM {\mathbf{M}}

\def \P {\mathcal{P}}

\def \arcm {{\mathrm{k}}}

\def \R {\mathbb{R}}

\def \der {s}
\def \derNode {j(\der)}

\def \setDERgi {\mathcal{S}_{\text{gi}}}
\def \setDERgf {\mathcal{S}_{\text{gf}}}
\def \setDERgfutil {\mathcal{S}_{\text{gf}}^{\text{utility}}}
\def \setDERgffacility {\mathcal{S}_{\text{gf}}^{\text{facility}}}
\def \setDERpq {\mathcal{S}_{\text{pq}}}
\def \setDERpqfixed {\mathcal{S}_{\text{pq}}^{\text{fixed}}}
\def \setDERpqvar {\mathcal{S}_{\text{pq}}^{\text{var}}}

\def \x {{x}}

\def \j {\mathbf{j}}

\def \state {\eta}
\def \opr {0}
\def \con {1}
\def \attack {a}
\def \optimalAttack {oa}
\def \defend {d}
\def \hardMin {hmin}
\def \hardMax {hmax}

\def \NN {{\mathrm{N}}}

\def \C {{\text{W}}}

\def \load {\text{LC}}
\def \shed {\text{LS}}

\def \Cshed {\C_{\shed}}
\def \Cload {\C_{\load}}
\def \Cmicro {\C_{\micro}}

\def \Clovr {\C_{\slovr}}
\def \Clofr {\C_{\slofr}}

\def \cost {{L}}
\def \costMicro {{L}_\microSmall}

\def \setPeriods {\mathcal{T}}

\def \noResponse {\text{\fontsize{5}{5}\selectfont AD}}
\def \micro {\text{\fontsize{5}{5}\selectfont MG}}
\def \maxMin {\text{\fontsize{6}{5}\selectfont Mm}}
\def \microSmall {\text{\fontsize{6}{5}\selectfont m}}
\def \noresponseSmall {{\text{\fontsize{7}{5}\selectfont nr}}}

\def \Second {\mathcal{D}_\arcm}
\def \ThirdMicro {{\Third}_\microSmall}
\def \SecondMicro {\Second^\microSmall}

\def \loss {\mathcal{L}}
\def \lossMaxmin {\mathcal{L}_{\maxMin}}
\def \lossMicro {\mathcal{L}_{\micro}}
\def \lossres {\mathcal{L}_{\text{\fontsize{8}{5}\selectfont res}}}
\newcommand{\lossNoResponse}{\loss_{\noResponse}}

\def \ltarget {\loss_\text{target}}
\def \resilience {\mathcal{R}}
\def \resilienceTarget {\resilience_\text{target}}
\def \resilienceMaxmin {\resilience_{\maxMin}}
\def \resilienceNoResponse {\resilience_{\noResponse}}
\def \resilienceMicro {\resilience_{\micro}}

\def \lcompleteShed {\loss_\mathrm{max}}

\def \y { y}

\def \i  {I}

\def \slovr {\text{VR}}
\def \slofr {\text{FR}}

\def \xx {X}

\def \l {\mathcal{L}}

\def \a {\delta}

\def \Xnpf {\mathcal{X}}

\def \ref {\text{ref}}

\def \second {{d}}

\def \third {u}

\newcommand{\Third}[1][]{{ \mathcal{U}^{#1}}}

\def \x {x}

\newcommand{\resistance}[1]{\mathbf{r}_{#1}} 
\newcommand{\reactance}[1]{\mathbf{x}_{#1}}

\def \f {\omega}

\def \ssum {\textstyle\sum}

\newcommand{\argmin}{\operatornamewithlimits{argmin}}

\usepackage{chngcntr}
\counterwithout{subsubsection}{subsection}

\titleformat{\subsubsection}[runin]
  {\normalfont\normalsize\itshape}{\thesubsubsection)}{2pt}{}[ - ]
\makeatletter
\@addtoreset{subsubsection}{subsection}
\makeatother

\titleformat{\paragraph}[runin]
{\normalfont\normalsize\itshape}{\thesubsubsection)}{2pt}{}[:]

\usepackage{tcolorbox}

\usepackage{mdframed}

\DeclareDocumentCommand{\mycommand}{ O{mydefault} m o o o }{%
	p:#2%
	\IfNoValueTF{#3}%
	{}%
	{p:#3}%
	\IfNoValueTF{#4}%
	{}%
	{\@ifmtarg{#4}{}{ p:#4}}%
	\IfNoValueTF{#5}%
	{}%
	{\@ifmtarg{#5}{}{ p:#5}}
	p:#1
}

\newcommand{\test}[3][o o o]{%
	\ifthenelse{\equal{#1}{}}{omitted}{given}%
}

\newcommand{\fmin}[1]{\bm{f^{min}_{#1}}} 
\newcommand{\fmax}[1]{\bm{f^{max}_{#1}}} 
\def \pre {o}
\def \post {c}

\def \edge {ij}
\def \period {t}

\def \scenarioIdx {s}
\def \node {i}
\def \ts {\period}
\def \nperiod {\mathrm{T}}

\newcommand{\myc}[3]{
	\def \firstString {{#1}}	
	\IfEqCase{#2}{
		{}{\firstString_{#3}}
		{\opr}{\firstString_{#3}^{\opr}}
		{\con}{\firstString_{#3}^{\con}}
		{\star}{\firstString_{#3}^{\star}}
		{\ts}{\firstString_{#3}^{\ts}}
		{\scenarioIdx}{\firstString_{#3}^{\scenarioIdx}}
		{\ts-1}{\firstString_{#3}^{\ts-1}}
		{0}{\firstString_{#3}^{0}}
		{\period}{\firstString_{#3}^{\period}}
		{\period-1}{\firstString_{#3}^{\period-1}}
		{\nperiod}{\firstString_{#3}^{\nperiod}}
		{n}{\firstString_{#3}}
		{l}{\widehat{\firstString}_{#3}}
		{u}{\widecheck{\firstString}_{#3}}
		{nr}{\firstString^{\noresponseSmall}_{#3}}
		{max}{\mathbf{\overline{\firstString}}_{#3}}
		{min}{\mathbf{\underline{\firstString}}_{#3}}
		{\hardMax}{\bm{\overline{\overline{\firstString}}_{#3}}}
		{\hardMin}{\bm{\underline{\underline{\firstString}}_{#3}}}	
		{constant}{\mathbf{\firstString}_{#3}}
		{\pre}{\firstString_{#3}^{\pre}}
		{\post}{{\firstString}_{#3}^{c}}
		{act}{\firstString_{#3}^{act}}
		{set}{\firstString_{#3}^{set}}
		{ref}{\mathbf{\firstString}_{#3}^{\text{ref}}}
		{r}{\firstString_{#3}^{r}}
		{nom}{\mathbf{\bm{\firstString}_{#3}^{nom}}}
		{stab}{\bm{\firstString_{#3}^{stab}}}
		{devmax}{\firstString_{#3}^{dev,max}}
		{reg}{\bm{\firstString_{#3}^{reg}}}
		{ev}{\firstString_{#3}^{ev}}
		{nev}{\firstString_{#3}^{nev}}
		{maxev}{\bm{\overline{\firstString}_{#3}^{ev}}}
		{maxnev}{\bm{\overline{\firstString}_{#3}^{nev}}}
		{\state}{\firstString_{#3}^{\state}}
		{\attack}{\firstString_{#3}^{\attack}}
		{\optimalAttack}{\firstString_{#3}^{\attack\star}}
		{\defend}{\firstString_{#3}^{\defend}}
	}
}

\newcommand{\fcnom}{\mathbf{f^{nom}}}

\newcommand{\uc}[2]{\myc{u}{#1}{#2}}

\newcommand{\ngc}[2]{\myc{\mathrm{G}}{#1}{#2}}

\newcommand{\kcc}[2]{\myc{kc}{#1}{#2}}	
\newcommand{\kgc}[2]{\myc{kg}{#1}{#2}}	
\newcommand{\klinec}[2]{\myc{kl}{#1}{#2}}	
\newcommand{\krc}[2]{\myc{kr}{#1}{#2}}

\newcommand{\nuc}[2]{\myc{\mathrm{v}}{#1}{#2}}
\newcommand{\nucref}[1]{\mathbf{v}^{\text{ref}}_{#1}}
\newcommand{\nucmax}[1]{\mathbf{\overline{v}}_{#1}}
\newcommand{\nucmin}[1]{\mathbf{\underline{v}}_{#1}}

\newcommand{\Pc}[2]{\myc{P}{#1}{#2}}
\newcommand{\Qc}[2]{\myc{Q}{#1}{#2}}

\newcommand{\nucc}[2]{\myc{vc}{#1}{#2}}
\newcommand{\nugc}[2]{\myc{vg}{#1}{#2}}
\newcommand{\fcc}[2]{\myc{fc}{#1}{#2}}
\newcommand{\fgc}[2]{\myc{fg}{#1}{#2}}
\newcommand{\pcc}[2]{\myc{pc}{#1}{#2}}
\newcommand{\qcc}[2]{\myc{qc}{#1}{#2}}

\newcommand{\ptc}[2]{\myc{p}{#1}{#2}}
\newcommand{\qtc}[2]{\myc{q}{#1}{#2}}

\newcommand{\Gnc}[1]{\myc{Gn}{}{#1}} 
\newcommand{\hnc}[1]{\myc{hn}{}{#1}}
\newcommand{\Gec}[1]{\myc{Ge}{}{#1}} 
\newcommand{\Hec}[1]{\myc{H}{}{#1}} 
\newcommand{\hec}[1]{\myc{he}{}{#1}}

\newcommand{\kpc}[1]{\myc{mp}{constant}{#1}}
\newcommand{\kqc}[1]{\myc{mq}{constant}{#1}}

\newcommand{\pegc}[2]{\myc{pe}{#1}{#2}}
\newcommand{\qegc}[2]{\myc{qe}{#1}{#2}}

\newcommand{\pgc}[2]{\myc{pg}{#1}{#2}}
\newcommand{\qgc}[2]{\myc{qg}{#1}{#2}}

\newcommand{\vdc}[2]{\myc{{\Delta \mathrm{v}}}{#1}{0}} 
\newcommand{\fdc}[2]{\myc{\Delta f}{#1}{0}} 
\newcommand{\prc}[2]{\myc{pr}{#1}{#2}}
\newcommand{\qrc}[2]{\myc{qr}{#1}{#2}}

\newcommand{\pnc}[2]{\myc{pn}{#1}{#2}}
\newcommand{\qnc}[2]{\myc{qn}{#1}{#2}}

\newcommand{\xc}[2]{\myc{\x}{#1}{#2}}
\newcommand{\Xc}[2]{\myc{\Xnpf}{#1}{#2}}
\newcommand{\Xcmicro}[2]{\myc{\Xnpf}{#1}{\microSmall}}
\newcommand{\Ycmicro}[2]{\myc{\mathcal{Y}}{#1}{\microSmall}}
\newcommand{\lcc}[2]{\myc{\beta}{#1}{#2}}

\newcommand{\fc}[2]{\myc{f}{#1}{#2}}

\newcommand{\transpose}[1]{{#1}^{\top}}

\newcommand{\mycc}[4]{
	\def \firstString {{#1}}	
	\IfEqCase{#2}{
		{n}{\firstString_{#3}^{#4}}
		{max}{\overline{\firstString}_{#3}^{#4}}
		{min}{\underline{\firstString}_{#3}^{#4}}
	}
}

\begin{document}
	\title{Resilience of Electricity Distribution Networks \\
		Part II: Leveraging Microgrids }
	\author{Devendra Shelar, Saurabh Amin, and Ian Hiskens
	\thanks{Manuscript submitted on April 10, 2019. This work was supported by NSF project ``FORCES'' (award $\#$: CNS-1239054), NSF CAREER (award $\#$: CNS-1453126), and NSF project \enquote{Modeling and Analysis of Load Ensembles} (award $\#$: ECCS-1810144).}
	
	\thanks{Devendra Shelar and Saurabh Amin are with Department of Civil and Environmental Engineering, Massachusetts Institute of Technology, 77 Massachusetts Avenue 1-241, Cambridge, MA 02139 USA (e-mail: \texttt{shelard,amins}@mit.edu, phone: 857-253-8964).}
	\thanks{I. A. Hiskens is with the Department of Electrical Engineering and Computer	Science, University of Michigan, Ann Arbor, MI 48109 USA (e-mail: hiskens@umich.edu).}
}
	\maketitle	

	\begin{abstract}
		Advances in microgrids powered by Distributed Energy Resources (DERs) make them an attractive response capability for improving the resilience of electricity distribution networks (DNs). This paper presents an approach to evaluate the value of implementing a timely response using microgrid operations and DER dispatch in the aftermath of a disruption event, involving strategic compromise of multiple DN components. Firstly, we extend the  resiliency assessment framework in \cite{part1} and develop a sequential (bilevel) model of attacker-operator interactions on a radial DN with one or more microgrids. Particularly, the operator response  includes microgrid operations under various islanding configurations (regimes), and single- or multi-master operation of DERs in providing grid-forming services as well as frequency and voltage regulation. Secondly, we introduce a restoration problem in which the operator gradually reconnects the disrupted components over multiple periods to restore the nominal performance of the DN. The first problem, formulated as a bilevel mixed-integer problem, is solved using a Benders decomposition method. The second problem, formulated as a multi-period mixed-integer problem, can be solved using a greedy algorithm. Our results illustrate the benefit of using microgrids in reducing the operator's losses, both immediately after the disruption event and during the restoration process. 
		
	\end{abstract}
	
	\section{Introduction}\label{sec:introduction2} 
	Modern electricity Distribution Networks (DNs) are prone to risks of service interruptions due to the failures of unreliable, and often insecure, cyber and physical components. Recent disruptions caused by natural disasters~\cite{microgridDisasterRecovery} and security attacks~\cite{ukraine,microgridSecurity} highlight the vulnerability of DNs to cyberphysical failures. 
	In this article, we investigate the use of microgrid technologies such as microgrid islanding and dispatch of Distributed Energy Resources (DERs)~\cite{microgridEmergencyControl,microgridReserveManagement,microgridSecurity} toward  improving DN resilience. 
	Historically, the idea of DER-powered microgrids as a response mechanism has been considered for responding to reliability failures~\cite{microgridReserveManagement,microgridSecurity,microgridsNikos}. Indeed, microgrids have been implemented to support the reliability targets of critical facilities such as hospitals, industrial plants, and military bases. However, their technological feasibility (and related operational aspects) in responding to security failures has received limited attention. We address this issue by building on our work in \cite{part1}, and focus on evaluating the effectiveness of DER-powered microgrids in limiting post-contingency losses after a disruption. 
	
	%
	
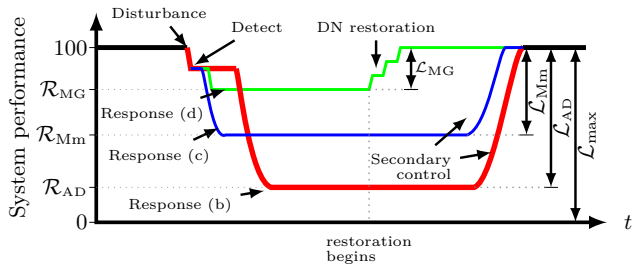
\begin{figure}[htbp!]
		\tikzstyle{mfs} = [font=\fontsize{9}{8}\selectfont,align=left]
		\tikzstyle{mrfs} = [mfs,rotate=90]
		\tikzstyle{fs} = [font=\fontsize{6}{6}\selectfont,align=left]
		\tikzstyle{ft} = [draw,line width=2pt]
		\tikzstyle{ftr} = [draw,line width=2.2pt,red]
		\tikzstyle{ftb} = [draw,line width=1.2pt, blue]
		\tikzstyle{ftg} = [draw,line width=1.2pt, green]
		\tikzstyle{ptr} = [-latex,line width=1pt]
		\tikzstyle{mk} = [gray,line width=0.6pt, dotted, opacity=0.8]
		\def \maxp {2.5}
		\def \dip {1}
		\def \rebound {2}
		\resizebox{8.5cm}{!}{\begin{tikzpicture}
			
			\node[](ytip) at (-0.5,3) {}; \node[](xtip) at (6.8,0) {};
			\draw[<->,>=latex,line width=1.5pt] (ytip.center) |- (xtip.center); 
			\node[mfs,right=0 of xtip]() {$t$}; 
			
			\draw[ft] (-0.5,2.5) -- (0.8,2.5) (5.6,2.5) -- (6.5,2.5);
			\draw[ftr] (0.8,2.5) -- (0.85,2.2) -- (1.5,2.2);
			\path[ftg] (0.85,2.2) -- (1.1,2.2) -- (1.15,1.9) -- (3.4,1.9) -- (3.45,2.1) -- (3.6,2.1) -- (3.65, 2.3) -- (3.8,2.3) -- (3.85,2.5) -- (5.6,2.5);
			\path[ftr] (1.5,2.2) .. controls  (1.55,2) and (1.7,0.6) .. (2,0.5) -- (4.9,0.5) .. controls (5.2,0.6)  and (5.45,2.45) .. (5.6,2.5);
			\path[ftb] (0.85,2.2) -- (1.00,2.2) .. controls (1.05,2.2) and (1.15,1.15) .. (1.35,1.25) -- (4.8,1.25) .. controls (5.1,1.35) and (5.25,2.45) .. (5.35,2.5) --(5.6,2.5);
			
			\node[fs] (dn) at (0.4,3) {Disturbance}; 
			\node[fs] (den) at (1.8,2.8) {Detect}; 
			\node[fs] (cas) at (0.7,0.25) {Response (b)};
			\node[fs] (trn) at (0.40,0.95) {Response (c)}; 
			\node[fs] (tcn) at (0.30,1.55) {Response (d)};

			\node[fs] (rn) at (3.5,2.8) {DN restoration}; 
			\node[fs,align=right] (svn) at (4,0.8) {Secondary\\control}; 
			\draw[ptr] (dn) -- (0.8,2.5);
			\draw[ptr] 	(rn) -- (3.5,2.2);
			\draw[ptr] 	(tcn) -- (1.35,1.9);
			\draw[ptr] 	(trn) -- (1.25,1.35);
			\draw[ptr] 	(den) -- (0.9,2.25);
			\draw[ptr] 	(svn) -- (5.15,1);
			\draw[ptr] 	(svn) -- (4.8,1.5);
			\draw[ptr] 	(cas) -- (1.9,0.45);
			\def \a {-1.2}; \def \b {-0.5}; \def \c {2.5}; \def \d {1.25}; \def \e {0.5}; \def \f {0}; \def \g {-0.8}; \def \h {-1.3}; \def \i {-1.8};
			\def \yredlow {0.5}; 
			\def \ybluelow {1.25}; 
			\def \ygreenlow {1.9};
			\def \axmicro {3.4}; 
			\def \axmaxmin {4.8}; 
			\def \axnoresponse {4.9};
			\def \axxmicro {4.1}; 
			\def \axxmaxmin {5.65}; 
			\def \axxnoresponse {6};  
			\def \axxmax {6.35}; 
			\def \ymax {2.5};
			\def \start {-0.5};
			\foreach \y/\l in {0/0,\yredlow/\resilienceNoResponse, \ybluelow/\resilienceMaxmin, \ygreenlow/\resilienceMicro, \ymax/100}{
				\draw (-0.6,\y) -- (-0.5,\y);
				\node[mfs,anchor=east] () at (-0.5,\y) {$\l$}; 
			}
			\foreach \x/\y/\z in {\start/\yredlow/2,\start/\ybluelow/1.2,\start/\ygreenlow/1.05}{
				\draw[gray,opacity=0.8,dotted] (\x,\y) -- (\z,\y);
			}
			
			\node[mrfs] (ylabel) at (-1.6,1.5) {System performance};
			
			\tikzstyle{ptre} = [ptr,<->,>=latex]

			\def \w {0.2}
			\foreach \x/\y/\l in {\axxmaxmin/\ybluelow/\lossMaxmin, \axxnoresponse/\yredlow/\lossNoResponse,\axxmax/0/\lcompleteShed} {
				\pgfmathsetmacro\xx{\x - \w/2}
				\pgfmathsetmacro\xxx{\x + \w/2}
				\draw [ptre] (\x,\ymax) -- (\x,\y);
				\draw [] (\xx,\y) -- (\xxx,\y);
				
				\pgfmathsetmacro\lx{\x + 0.165}
				\pgfmathsetmacro\ly{(\y + \ymax)/2}
				\node[mrfs] (lnr) at (\lx,\ly)  {$\l$};
			}
			\foreach \xx/\y/\xxx in {\axmicro/\ygreenlow/\axxmicro, \axmaxmin/\ybluelow/\axxmaxmin, \axnoresponse/\yredlow/\axxnoresponse} {
				\draw [mk] (\xx,\y) -- (\xxx,\y);
			}
		\draw [mk] (3.4,1.9) -- (3.4,0);
		\node[fs] (restoration) at (3.4,-0.4)  {restoration\\begins};
			
			\draw[] (4.1,\ygreenlow) -- (3.9,\ygreenlow);
			\draw[ptre] (4,\c) -- (4,1.9);
			\node[fs] (lmg) at (4.35,2.2)  {$\lossMicro$};
			\end{tikzpicture}}
		\caption{Performance under various response capabilities.} \label{fig:resilienceDefinition2}
	\end{figure}

	We model the  sequential interaction between a DN operator and an external adversary as follows~\cite{part1}: 
	\begin{alignat}{10}\label{eq:maxMinGeneric2}\tag{P1}
	\lossMaxmin &\coloneqq\quad && \max_{\second\in\Second} &\min_{\third\in\Third[](\second)} &&\ \cost\left(\third,\xc{}{}\right) \qquad&\text{s.t.} ~\quad&& \xc{}{} \in \Xc{}{}\left(\third\right),
	\end{alignat}
	where $\second\in\Second$ denotes an attacker strategy, $\third\in\Third(\second)$ an operator response strategy, $\xc{}{}\in \Xc{}{}$ the network state, and $\cost$ the composite loss function, the details of which are presented in \cref{sec:bilevel}.  In \cite{part1}, we argued that cyberphysical disruptions to DNs can lead to operating bound violations and cause uncontrolled or forced disconnects of DN components. Specifically, we modeled the impact of attacker-induced disconnects of DN components as supply-demand disturbances, and the impact of TN-side disturbances as voltage deviations at the substation node. Then, we considered preemptive load control and component disconnects as operator response actions for the generic setting when the attacker's (resp. operator's) goal is to maximize (resp. minimize) the post-contingency losses. We introduced $\resilienceMaxmin \coloneqq 100\left(1-\lossMaxmin/\lcompleteShed\right)$ as a resilience metric of the DN, where $\lcompleteShed$ (chosen for sake of normalization) denotes the operator loss when all DN components are disconnected; see  \Cref{fig:resilienceDefinition2}. 
Finally, we evaluated the value of optimal response as the total reduction in post-contingency losses relative to the case of autonomous (local protection driven) disconnections, i.e. $\resilienceMaxmin-\resilienceNoResponse$, where $\resilienceNoResponse = 100\left(1-\lossNoResponse/\lcompleteShed\right)$.

	In this article, we consider another bilevel formulation: 
	\begin{alignat}{10}\label{eq:maxMinMicroGeneric}\tag{P2}
	\lossMicro &\coloneqq\quad && \max_{\second\in\SecondMicro} \ &\min_{\third\in\ThirdMicro(\second)} &&\ \costMicro\left(\third,\xc{}{}\right) \ \ &\text{s.t.} \ && \xc{}{} \in \Xcmicro{}{}\left(\third\right),
	\end{alignat}
	where the network model $\Xcmicro{}{}$, and the loss function $\costMicro$ are extended to capture  microgrid operations (\cref{sec:microgridModel}) and DER dispatch and regulation aspects (\cref{sec:derModel});
	 the set of attacker strategies $\SecondMicro$ and the set of operator strategies $\ThirdMicro$ are 
	 also modified to capture attacker-operator interactions for DNs with DER-powered microgrids~(\cref{sec:bilevel}). The maximin value of \eqref{eq:maxMinMicroGeneric}, $\lossMicro$, denotes the worst-case post-contingency loss incurred by the operator for the given microgrid and DER capabilities; see  \Cref{fig:resilienceDefinition2}. Then,  $\resilienceMicro \coloneqq 100\left(1-\lossMicro/\lcompleteShed\right)$ can be viewed as a resilience metric of the DN under microgrid-enabled operator response.  
	\ifArxivVersion
		Furthermore, the relative value of timely microgrid response (or equivalently, the improvement in DN resilience due to microgrids) can be evaluated as $(\resilienceMicro - \resilienceMaxmin)$. We posit that advances in DER-enabled microgrids and emergency control operations at the substation level can be leveraged to implement timely resiliency-improving response actions (less than a few seconds after a disturbance event). 
	\fi
	
	In \cite{part1}, we considered three operator response models: 
	\begin{itemize}
		\item[(a)] Remote control by the control center during normal conditions; 
		\item[(b)]  Autonomous (protection driven) disconnection of individual components (tripping of DGs or loads under nodal violations in operating conditions); and 
		\item[(c)] Emergency control by a secure Substation Automation (SA) system. 
	\end{itemize}	
	As in \cite{part1}, it is assumed that the emergency control actions (c) subsume the autonomous actions (b) by making decisions that are coordinated across the SA system. Hence, (b) and (c) are never simultaneously active. Furthermore, in this paper, we consider the following extension of (c):
	\begin{itemize}
		\item[(d)] Emergency control by the SA involving microgrid islanding and DER dispatch. 
	\end{itemize} 
	Analogous to \cite{part1}, we consider that the SA system can detect the disrupted components from  changes in measurements of net nodal consumption. By using knowledge of the attack the SA can compute and implement the operator response in a timely manner. For our purposes, response (d) is an optimal second-stage response in \eqref{eq:maxMinMicroGeneric}. Our analysis relies on the premise that such response can be implemented via modern SA systems during disruptions. Indeed, the continued improvements in SA system disturbance detection and control capabilities can further assist in restoration operations.   
	
	The resiliency of a system also reflects how quickly it can rebound to its nominal state after a disruption~\cite{part1,resilienceDefnNIAC}. Microgrids can provide partial demand satisfaction during the system restoration process,  especially during the time when the DN is fully disconnected from the TN. We consider an admittedly simple, but practically relevant, multi-period DN restoration problem in which the disrupted DN components are gradually restored over several periods; refer to \enquote{DN restoration} in \Cref{fig:resilienceDefinition2}. Our goal in this problem is to compute an  operator strategy in each time period (roughly, on the order of a few minutes). Such a strategy is comprised of reconnecting disrupted components, modifying the microgrid islanding configuration, and dispatching the DERs within individual microgrids. 
	
	Our modeling approach addresses some key issues regarding microgrid and DER operations. In particular, we allow for the formation of one or more microgrid islands in radial DNs. When all the microgrids are connected to the transmission network (TN), the DN is operating in the \emph{grid-connected} regime. If none of the microgrids are connected to the TN, then the DN is operating in the \emph{fully-islanded} regime. In our approach, the DN can also operate in a \emph{partially-islanded} regime, in which some of the microgrids are connected to the TN while other microgrids are not. In both partially- and fully- islanded regimes, each microgrid can operate as an isolated microgrid or as part of a larger microgrid. To model power flows in each of the microgrids, we introduce a natural extension of the LinDistFlow equations. The resulting network model captures DN operations in all the above-mentioned regimes (\cref{sec:microgridModel}). We limit attention to linear power flows mainly for the ease of exposition.  
	
	Importantly, we consider the parallel operation of multiple DERs for the provision of  \emph{grid-forming} services, which involve providing voltage and frequency references, as well as maintaining voltage and frequency within operating bounds (i.e. \emph{regulation} services). When a microgrid is connected to the TN, the bulk generators provide the grid-forming services. However, when a microgrid is disconnected from the TN, then at least one  DER within that microgrid must provide the grid-forming services~\cite{microgridEmergencyControl}. 
	Depending on the number of grid-forming DERs within a microgrid, one can consider two modes of DER  operation under islanded regimes, namely Single-Master Operation (with a single grid-forming DER) and Multi-Master Operation (with more than one grid-forming DER)~\cite{microgridEmergencyControl}. Our model is sufficiently flexible to capture both the single- and multi-master modes of DER operation. 
	In addition to voltage regulation, we also consider frequency regulation, which becomes important for microgrids due to the low inertia of the DERs. By using the appropriate droop control equations, we capture both frequency and voltage regulation aspects resulting from multiple DERs operating in parallel within a microgrid~(\cref{sec:derModel}). 
	
Our main contributions are as follows:
\begin{itemize}[label=$(\star)$]
	\item We capture the different microgrid regimes as well as DER operating modes  by developing a new mixed-integer linear network model. This modeling approach enables us to formulate \eqref{eq:maxMinMicroGeneric} as a Bilevel Mixed-Integer Problem (BiMIP). In~\cite{part1}, we showed that \eqref{eq:maxMinGeneric2} is also a BiMIP, and can be solved using a Benders Decomposition (BD) algorithm. In \cref{sec:bilevel}, we show  that this algorithm can be applied to the extended formulation \eqref{eq:maxMinMicroGeneric}. 
	\item Our network model is also well-suited for formulating a DN restoration problem as a multi-period Mixed-Integer Problem (MIP). In our restoration problem, the network state in any period only depends on the operator response actions in that period, and the network state in the previous period. We exploit this feature and propose a greedy heuristic that seeks to reconnect the disrupted components in each period such that the post-contingency losses for that period are minimized (\cref{sec:restoration}). 
\end{itemize}	

	\section{Multi-Microgrid DN  model}\label{sec:microgridModel}
	
	In this section, we develop a model of a radial DN with one or more microgrids. 
	This network model extends the LinDistFlow model~\cite{lindist} to multi-microgrid settings. 
	
	We distinguish between two operating stages $\pre$ and $\post$, which denote the pre- and post- contingency stage, respectively. The network is initially in $\pre$ stage, and after the disturbance event enters the $\post$ stage; see~\cref{sec:bilevel} for details on the disturbance model. 
	Let $\state\in\{\pre,\post\}$ denote the operating stage of the network. We define the network state as $\xc{\state}{} \coloneqq  \transpose{\left(\ptc{\state}{}, \qtc{\state}{}, \Pc{\state}{},\Qc{\state}{}, \nuc{\state}{}, \fc{\state}{}\right)}$, where each of these entries are themselves vectors of appropriate dimensions, and are described in \cref{tab:notationsTable2}. 
	
			\begin{table}[htbp!]
		\renewcommand{\baselinestretch}{1}
		\caption{Table of Notation.} 
		\label{tab:notationsTable2}
		\centering
		\def \xx {0.5cm}
		\def \xxx {5.0cm}
		\def \xxxx {7.0cm}
		\def \xxxxx {1cm}
		\def \xxxxxx {10cm}
		\def \xxxxxxx {6.5cm}

		\resizebox{8.95cm}{!}{
			\begin{tabular}{p{2.5cm}p{8.8cm}}
				\multicolumn{2}{l}{\textbf{DN  parameters}} \\
				\ifShortTable
				$\N$ & set of nodes in DN \\
				$\E$ & set of edges in DN \\
				$0$ & substation node label\\
				\fi	
				$\setMicrogrid\subseteq \E$ & set of microgrid connecting lines\\
				$\N[i]\subseteq \N$ & nodes belonging to $i^{th}$ microgrid\\
				$\setMicrogrid_{i}\subseteq \setMicrogrid$ & set of lines which if open isolate the $i^{th}$ microgrid\\
				\ifShortTable
				$\j$ &  complex square root of -1, $\j = \sqrt{ -1}$ \\
				$\nuc{nom}{}$ & nominal squared voltage magnitude (1 pu)\\
				\fi
				$\fcnom$ & nominal system frequency (1 pu)\\	
				$\P_i \subseteq \E$ & lines on the path between node $i$ and  substation node\\	
				\multicolumn{2}{l}{\textbf{DER categories}} \\				
				$\setDER$ & set of DERs\\
				$\setDERgf\subseteq\setDER$ & set of grid-forming DERs\\
				$\setDER\subseteq \setDERgf$ & set of PQ Inverter (PQI)-controlled DERs\\
				$\setDERpqfixed \subseteq \setDERpq$ & set of PQI-controlled DERs with fixed setpoints\\
				$\setDERpqvar \subset \setDERpq$ & set of PQI-controlled DERs with controllable setpoints\\
				$\setDERgfutil \subseteq \setDERgf$ & set of utility-owned grid-forming DERs\\
				$\setDERgffacility \subseteq \setDERgf$ & set of facility level microgrid-specific grid-forming DERs \\
				$\setDERgi = \setDERgf\bigcup\setDERpqvar$ & set of grid-interactive DERs\\
				\multicolumn{2}{l}{\textbf{Nodal quantities of node $i \in \N $}} \\
				\ifShortTable
				$\nuc{}{i}$ & squared voltage magnitude at node $i$\\
				\fi
				$\fc{}{i}$ & system frequency measured at node $i$ \\		
				\ifShortTable
				$\nucc{min}{i},\nucc{max}{i}$ & lower, upper voltage bounds for load $i$\\
				$\nugc{min}{i},\nugc{max}{i}$ & lower, upper voltage bounds for DG $i$\\
				\fi
				$\fcc{min}{i},\fcc{max}{i}$ & lower, upper frequency bounds for load $i$\\
				$\fgc{min}{i},\fgc{max}{i}$ & lower, upper frequency bounds for DG $i$\\
				\ifShortTable
				$\pcc{max}{i} +\j\qcc{max}{i}$ & nominal demand at node $i$  \\
				$\pcc{}{i} + \j \qcc{}{i}$ & actual power consumed at node $i$  \\
				$\kcc{}{i}\in\{0,1\}$ & 0 if load $i$ is connected to DN; 1 otherwise \\		
				$\lcc{}{i}$ & fraction of demand satisfied at node $i$\\
				$\lcc{min}{i}$ & lower bound of  load control parameter $\lcc{}{i}$\\
				$\pgc{max}{i} +\j\qgc{max}{i}$ & nominal generation of DG  $i\in\setDERpqfixed$ \\
				$\pgc{}{i} + \j\qgc{}{i}$ & actual power generated by DG $i\in\setDERpqfixed$ \\
				$\kgc{}{i}\in\{0,1\}$ & 0 if DG $i\in\setDERpqfixed$ is connected to DN; 1 otherwise \\
				\fi
				\multicolumn{2}{l}{\textbf{Quantities of DER $\der \in \setDER$}} \\
				$\derNode$ & the DN node where the DER $\der\in\setDERgf$ is located\\
				$J(S) \subseteq\N$ & the set of DN nodes where the DERs in the set $S\subseteq\setDER$ are located\\ 
				$\pnc{min}{\der}, \pnc{max}{\der}$ & lower, upper active power bounds of microsource $\der\in\setDERgf$\\
				$\qnc{min}{\der}, \qnc{max}{\der}$ & lower, upper reactive power bounds of microsource $\der\in\setDERgf$\\
				$\pnc{}{\der}+\j\qnc{}{\der}$ & total power supplied by microsource $\der\in\setDERgf$ \\
				$\pegc{min}{\der}, \pegc{max}{\der}$ & lower, upper active power bounds of storage $\der\in\setDERgf$\\
				$\qegc{min}{\der}, \qegc{max}{\der}$ & lower, upper reactive power bounds of storage $\der\in\setDERgf$\\
				$\pegc{}{\der}+\j\qegc{}{\der}$ & total power supplied by storage $\der\in\setDERgf$ \\
				$\krc{}{\der}\in\{0,1\}$ & 1 if DER $\der\in\setDERgf$ contributes to grid-forming services\\
				$\prc{}{\der}+\j\qrc{}{\der}$ & total power supplied by DER $\der\in\setDERgf$ \\
				$\prc{ref}{\der}, \qrc{ref}{\der}$ & active, reactive power references of DER $\der\in\setDERgf$\\
				$\fc{ref}{\der}, \nucref{\der}$ & frequence, voltage references of DER $\der\in\setDERgf$\\
				\multicolumn{2}{l}{\textbf{Parameters of edge $(i,j) \in \E $}}  \\
				$\klinec{}{\edge}\in\{0,1\}$ & 1 if $(i,j)$ is switched open; 0 otherwise  \\
				\ifShortTable
				$\Pc{}{\edge} +\j\Qc{}{\edge} $ & power flowing from node $i$ to node $j$ \\
				$\resistance{\edge}, \reactance{\edge}$ & resistance and reactance of line $(i,j) \in \E$  \\
				\fi
			\end{tabular}
		}	
	\end{table}


	
\begin{figure}[htbp!]
	\centering
	\tikzset{every node/.append ={font size=tiny}} 
	\tikzstyle{dnnode}=[draw,circle, minimum size=0.5pt, inner sep = 2]
	\tikzstyle{dnedge}=[-, line width=1pt]
	\tikzstyle{dernode}=[circle, fill=blue, minimum size=0.5pt, inner sep = 2]
	\tikzstyle{blackoutnode}=[circle, fill=black, minimum size=0.5pt, inner sep = 2]
	\tikzstyle{failededge}=[-, densely dotted]
	\def \drawgrid {\draw[step=1,gray, ultra thin, draw opacity = 0.5] (0,0) grid (3,4);}
	\def \drawSubstation {\draw[-, line width = 2pt] (0.8,4.05) -- (2.2,4.05)  node [midway,above] {};}
	\def \drawZero {\node (0) at (1.5,3.85) {};}
	\def \drawBus at (#1,#2); {\draw [line width=2.5pt,-] (#1,#2-0.6) -- (#1,#2+0.6);}
	\def \drawBuss at (#1); {\draw let \p1=(#1) in [draw,line width=2.5pt,-] (\x1,\y1+20) -- (\x1,\y1-20);}
	\def \drawLines at (#1,#2,#3,#4); {\draw [line width=1.5pt,-] (#1,#2) -- (#3,#4);}
	\def \drawArrows at (#1,#2,#3,#4); {\draw [line width=1.5pt,->, >=latex] (#1,#2) -- (#3,#4);}
	\def \switch (#1,#2); {
		\path (#1) ++(40:1) node (#1a) {};
		\path (#1) ++(0:0.8) node (#1b) {};
		\path (#1) ++(40:0.8) node (#1c) {};
		\draw[#2] (#1.center) -- (#1a.center); 
		\draw[->, line width=0.5pt,>=latex,#2] (#1b.center) to [out=90,in=-60] (#1c.center);
	}
	\def \rotateSwitch (#1,#2,#3,#4); {
		\path (#1) ++(#4*40+#4*#3:1) node (#1a) {};
		\path (#1) ++(#4*#3:0.8) node (#1b) {};
		\path (#1) ++(#4*40+#4*#3:0.8) node (#1c) {};
		\draw[#2] (#1.center) -- (#1a.center); 
		\draw[->, line width=0.5pt,>=latex,#2] (#1b.center) to [out=#4*90+#4*#3,in=-60*#4+#4*#3] (#1c.center);
	}
	\def \dis {1.8}
	\def \drawOval (#1,#2); {\node (n12oval) at ($(#1)!0.5!(#2)$) {};	
		\fill[gray,opacity=0.4] (n12oval)[xshift=\dis*0.5] circle [x radius=\dis*1.2, y radius=\dis/2*1.2, xshift=\dis*0.5];	} 	
	\resizebox{8.5cm}{!}{
		\begin{tikzpicture}[scale=0.6]
		\def \ly {1.5} 
		
		\node (substation) at (0,0) {}; 
		\node[right=\dis of substation](node1) {};
		\node[right=\dis of node1](node2) {};
		\node[above right=\dis*5/13 and \dis*12/13 of node2](node3) {};
		\node[below right=\dis*5/13 and \dis*12/13 of node2](node5) {};
		\node[right=\dis of node3](node4) {};
		\node[right=\dis of node5](node6) {};
		
		\foreach \xx [count=\i] in {substation,node1,node2,node3,node4,node5,node6}{
			\pgfmathsetmacro\result{\i - 1}
			\drawBuss at (\xx);
			\node[below= 0.4 of \xx] () {\pgfmathprintnumber{\result}};
		}
		
		\foreach \x/\y in {substation/node1,node1/node2,node3/node2,node5/node2,node4/node3,node6/node5} 
		\draw [line width=1.5pt,-]  (\x.center) -- (\y.center); 
		
		\drawOval (node1,node2); \drawOval (node3,node4); \drawOval (node5,node6);

		\node[right= 0.5 of substation](dnswitch) {}; \switch(dnswitch,black);		
		\node (n23mid) at ($(node2)!0.4!(node3)$) {};	\rotateSwitch(n23mid,black,25,1);	
		\node (n25mid) at ($(node2)!0.4!(node5)$) {};	\rotateSwitch(n25mid,black,25,-1);	
		
		\scalefont{0.7}
		\node[above=0.4 of node1,xshift=30] {$\N[1]=\{1,2\}$}; 
		\node[below=0.45 of node3,xshift=30] {$\N[2]=\{3,4\}$}; 
		\node[below=0.45 of node5,xshift=30] {$\N[3]=\{5,6\}$}; 
		\node[align=left,anchor=west](graph) at (0,-2.5) {$\setMicrogrid_1 = \{(0,1), (2,3), (2,5)\}$\\
			$\setMicrogrid_2 = \{(2,3)\}$\\
			$\P_5 = \{(0,1), (1,2), (2,5)\}$};
		\node[align=left,anchor=north] at (1.0,1.8) {Substation\\node}; 
		\node[align=left] at (8.1,1.8) {Connecting\\Line}; 
		
		\end{tikzpicture}
	}
	\caption{Multi-microgrid DN model. }
	\label{fig:systemStateDistribution2}
	
	\end{figure}
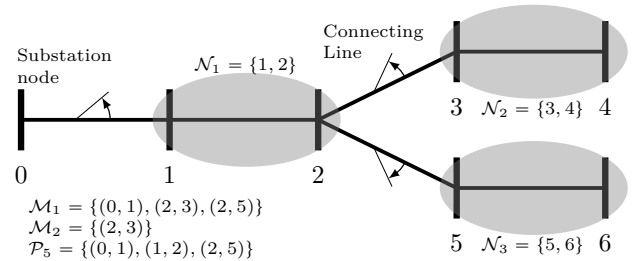

	In our DN model, we consider a radial network consisting of one or more microgrids. We refer to a distribution line $(i,j)\in\setMicrogrid\subseteq\E$ as a microgrid \emph{connecting line}  if it connects a microgrid to the TN or to other microgrids; see~\cref{fig:systemStateDistribution2}. Here $\setMicrogrid$ denotes a given fixed set of connecting lines. For a connecting line $(i,j) \in \setMicrogrid$, we use $\klinec{\state}{ij} = 0$  (resp. $\klinec{\state}{ij} = 1$) to indicate that it is in the closed (resp. open) state. The DN operating regimes are defined according to the states of the connecting lines: 
	\begin{itemize}[label=-]
		\item \emph{Grid-connected regime} when all connecting lines are closed (i.e $\klinec{\state}{ij} = 0 \quad \forall \ (i,j)\in \setMicrogrid$), 
		\item \emph{Fully-islanded regime} when all connecting lines are open (i.e. $\klinec{\state}{ij} = 1 \quad \forall \ (i,j)\in \setMicrogrid$), or 
		\item \emph{Partially-islanded regime} when there exists at least two connecting lines such that one of them is closed and the other is open (i.e. $\exists\ (i,j), (m,n) \in \setMicrogrid$ such that $\klinec{\state}{ij} = 0 \text{ and } \klinec{\state}{mn} = 1$).
	\end{itemize} 
	

	Let $\{ \N[1], \cdots, \N[\abs{\setMicrogrid}] \}$ denote the set of disjoint microgrid subnetworks of the DN, where each $\N[i]$ for $i \in \{1,\cdots,\abs{\setMicrogrid}\}$ denotes a  connected subnetwork when all connecting lines are open, i.e. $\klinec{\state}{mn} = 1$ for all $(m,n) \in \M$. For each subnetwork $\N[i]$, let $\setMicrogrid_{i}\subseteq\setMicrogrid$ denote the set of connecting lines which need to be open for $\N[i]$ to be completely isolated (i.e. autonomously operating). A \emph{microgrid island} is formed when an individual microgrid or a connected subnetwork of more than one microgrid no longer receives power supply from the TN. Also, let $\P_i$ denote the set of lines along the path connecting node $i$ to the substation node. For example, in \cref{fig:systemStateDistribution2} the set of connecting lines for the subnetwork $\N[1] = \{1,2\}$ is $\setMicrogrid_1 = \{(0,1), (2,3), (2,5)\}$. For this example, $\P_5 = \{(0,1),(1,2),(2,5)\}$. Also, if $\klinec{\eta}{01} = 1$, $\klinec{\eta}{23} = 0$, and $\klinec{\eta}{25} = 1$, then the microgrid $\N[3]$ is operating as an isolated island, whereas microgrids $\N[1]$ and $\N[2]$ are operating together as part of one larger microgrid island. 

\begin{remark}
	The smaller microgrids are typically used for supplying power to a critical facility (e.g. hospital, university, prison). In our model, these microgrids can be leveraged to supply power to the DN during emergency conditions (fully- or partially-islanded regimes).
\end{remark}

	Now, we describe the constraints related to the power flows, nodal frequencies and load connectivity in microgrids. Unless explicitly stated, the following constraints are valid for either operating stage $\state \in\{\pre,\post\}$. 

	\begin{enumerate}[wide, labelwidth=!, labelindent=10pt]
		\item \textit{Power flows:} 
		A connecting line permits power flow through it if and only if it is \textit{closed}.  We model this constraint as follows: 
		\begin{subequations}\label{eq:islandingCapacity}
			\begin{alignat}{8}
			\label{eq:islandingRealCap}	\abs{\Pc{\state}{\edge}} & \le \left(1-\klinec{\state}{\edge}\right)\bigM \qquad && \forall\ (i,j)\in\M\\
			\label{eq:islandingReacCap}	\abs{\Qc{\state}{\edge}} & \le \left(1-\klinec{\state}{\edge}\right)\bigM \qquad && \forall\ (i,j)\in\M,
			\end{alignat}
		\end{subequations}
		where $\bigM$ is a large constant. 
		This typical modeling trick to use a constraint of the type $\abs{a-c}\le y\bigM$ where $y\in\{0,1\}$,  enforces an equality $a = c$ only when $y=0$; otherwise the equality is not binding. We use this trick repeatedly to model various other constraints of a similar type. 
		
		\item \textit{Voltage drop:} 
		The voltage drop along a non-connecting line $(i,j) \not\in \setMicrogrid$ is given by the standard voltage drop equation of the LinDistFlow model~\cite{lindist}: 
		\begin{equation}\label{eq:voltageLindist}	
		\nuc{\state}{j} = \nuc{\state}{i} - 2\resistance{\edge}\Pc{\state}{\edge} - 2\reactance{\edge}\Qc{\state}{\edge}\quad \forall\ (i,j)\in\E\backslash\setMicrogrid. 
		\end{equation}
		However, for a connecting line, the voltage drop constraint is active only if it is closed, and is inactive, otherwise, i.e. 
		\begin{equation}\small	\label{eq:islandVoltConnected}	
		\begin{split}
		\abs{\nuc{\state}{j} - \left(\nuc{\state}{i} - 2\resistance{\edge}\Pc{\state}{\edge} - 2\reactance{\edge}\Qc{\state}{\edge}\right)} \le \klinec{\state}{\edge} \bigM \ \ \ 
		\forall \ (i,j) \in\M. 
		\end{split}
		\end{equation}
		
		\item \textit{Nodal frequencies:} 
	In islanded regimes, the DER(s) must provide grid-forming and regulation services~\cite{islandingControlStrategies,microgridEmergencyControl}. Moreover, a microgrid island can have multiple DERs operating in parallel. We assume that DERs can rapidly synchronize their frequencies to a common value with the help of power electronics~\cite{microgridEmergencyControl}. This value can be regarded as the island's frequency. To model that the nodal frequencies within a microgrid island are identical in steady state, we can write:
	\begin{equation*}
		\fc{\state}{i} = \fc{\state}{j}\quad \forall \ i,j \in\N[k] \text{ and } \forall\ k=1,\cdots,\abs{\setMicrogrid},
	\end{equation*}
	which is equivalent to writing: 
		\begin{equation}\label{eq:islandFrequencyConnected}	
		\fc{\state}{i} = \fc{\state}{j} \qquad \forall\ (i,j)\in\E\backslash\setMicrogrid, 
		\end{equation}
		because if a line $(i,j)$ is not a connecting line, i.e. $(i,j)\in\E\backslash\setMicrogrid$, then nodes $i$ and $j$ both belong to the same microgrid.
	Generically, frequency of every microgrid island can be different from the frequency of the TN-connected substation node. Moreover, the frequencies of any two microgrid islands that are not connected to each other can also be different. We model this constraint as follows:
		\begin{equation}\label{eq:islandFrequencyNotConnected}	
		\abs{\fc{\state}{i} - \fc{\state}{j}} \le \klinec{\state}{\edge} \bigM \qquad\forall\ (i,j)\in\setMicrogrid. 
		\end{equation}
		Finally, we model the constraint that the load gets disconnected (i.e. $\kcc{\state}{i} = 1$) when the nodal frequency violates the safe operating bounds:
		\begin{align}\label{eq:frequencyDisconnectLoad}
		\hspace{-0.4cm}		
		\begin{aligned}
		\kcc{\state}{i} &\ge \fcc{min}{i} - \fc{\state}{i}, \quad && 
		\kcc{\state}{i} \ge \fc{\state}{i} - \fcc{max}{i} \quad && \forall\ i\in \N.
		\end{aligned}
		\end{align}
		%
		%
		%
		
	\end{enumerate}
	
	
	
	
	\section{Distributed Energy Resources (DERs)}
	\label{sec:derModel}
\begin{figure}[htbp!]
			\tikzstyle{mnode} = [draw,align=left,minimum width=4cm]
			\tikzstyle{mnode2} = [mnode,minimum width=1.8cm]
			\tikzstyle{mcon} = [-,line width=2,>=latex]
			\tikzstyle{mcon2} = [mcon,->]
			\def \sep {0.25}
			\def \bsep {1}
			\def \lsepone {-1.2}
			\def \lseptwo {-1.5}
			\def \lsepthree {-1.3}
			\resizebox*{8.7cm}{!}{
				\begin{tikzpicture}
				[every node/.style={rounded corners},
				every edge/.style={black,-,thick,draw,line width = 4}]
				
				\node[mnode] (fullSet) {Set of DERs ($\setDER$)};				
				\node[mnode,below right =\bsep and \lsepone of fullSet] (pqi) {PQI-controlled ($\setDERpq$)};
				\node[mnode,below left =\bsep and \lsepone of fullSet] (gridForming) {Grid-forming ($\setDERgf$)};
				
				\node[mnode2,below left =\bsep and \lseptwo of gridForming] (gfutility) {($\setDERgfutil$)};
				\node[mnode2,below right =\bsep and \lseptwo of gridForming] (gffacility) {($\setDERgffacility$)};
				
				\node[mnode2,below left =\bsep and \lsepthree of pqi] (pqvar) { ($\setDERpqvar$)};
												
				\node[mnode2,below right =\bsep and \lsepthree of pqi, red, line width=2] (pqfixed) {($\setDERpqfixed$)};

				\node[below=0.2*\bsep of fullSet](fullSetBelow){};
				\node[below=0.2*\bsep of gridForming](gridFormingBelow){};
				\node[below=0.2*\bsep of pqi](pqiBelow){};
				
				\foreach \x/\y in {fullSet/fullSetBelow,gridForming/gridFormingBelow,pqi/pqiBelow}
					\draw [mcon] (\x.south) to (\y.center);
				
				\foreach \x/\y in {fullSetBelow/gridForming, fullSetBelow/pqi, gridFormingBelow/gfutility, gridFormingBelow/gffacility, pqiBelow/pqfixed, pqiBelow/pqvar}
					\draw [mcon2] (\x.center) -| (\y.north);
					
				\node[align=center,below =0.2 of gffacility] (gilabel){Grid-interactive ($\setDERgi$)};
				\node[mnode2,dashed,fit = {(gffacility) (gfutility) (gilabel) (pqvar)}, inner sep=7] (gnset) {};
				
\node[align=center,below =0.2 of pqfixed] (gnlabel){Grid-noninteractive}	{};
				\node[mnode2,dashed,fit = {(pqfixed) (gnlabel)}, inner sep=7] (gnset) {};
				\end{tikzpicture}
				
			}
			\vspace{0pt}
			\caption{Basic taxonomy of DERs~\cite{microgridsManagement}.}\label{fig:derClassification}
		\end{figure}
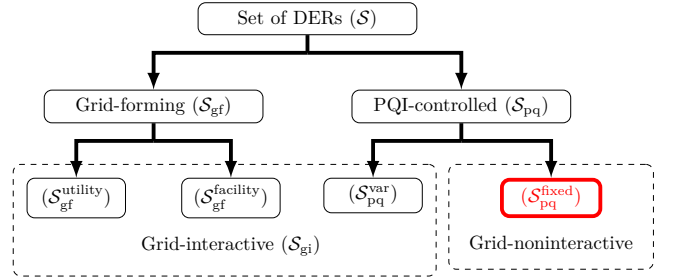

	We now introduce a generic taxonomy of DERs that is relevant to microgrid operations (see~\cite{microgridsManagement}) and a model which captures both single- and multi-master operating modes of DERs. 
	Refer to \cref{fig:derClassification} for DER categories and \cref{fig:derClassificationTab}  for a comparison of their capabilities. 
%
		
		\begin{table}[htbp!]
			\def\arraystretch{1.5}
		\def \mw {2.55cm}
		\def \mww {4.5cm}
		\def \mwww {4.5cm}
		\centering
		\setlength\extrarowheight{2pt}
		\subfloat[Grid-interactive vs. grid-noninteractive DERs.]{\label{tab:interactive}
			\resizebox{8.7cm}{!}{
				\begin{tabular}{x{\mw}|x{\mww}|x{\mwww}}
					\hline
					Attribute	& Grid-noninteractive ($\setDERpqfixed$) & Grid-interactive ($\setDERgi = \setDERgf\bigcup\setDERpqvar$) \\
					\hline\hline
					\parbox[c]{1.6cm}{Power \\output} & fixed & \parbox[l]{3cm}{variable/responsive to grid conditions}\\
					\hline
					\parbox{\mw}{Controllable \\ by response (a)} & \parbox[c]{4.1cm}{\hspace{1.7cm} yes\\ (remote set-point control)} & \parbox{\mwww}{no (can act as zero output source while being connected)}\\
					\hline
					\parbox{\mw}{Controllable \\ by response (b)} & \parbox[c]{3.2cm}{yes (due to operating bound violations)} & \parbox{\mwww}{no (LVRT \& LFRT available)}\\
					\hline
					\parbox{\mw}{Controllable by \\ response (c)} & yes & N/A\\
					\hline
					\parbox{\mw}{Controllable by \\ response (d)} & yes & yes\\
					\hline
		\end{tabular}}}\\
		\vspace{0pt}
		%
		\subfloat[PQI-controlled vs. grid-forming DERs.]{\label{tab:grid-forming}
			\resizebox{8.7cm}{!}{
				\begin{tabular}{x{\mw}|x{\mww}|x{\mwww}}
					\hline
					Attribute	& PQI-controlled DERs with variable setpoints ($\setDERpqvar$)& Grid-forming ($\setDERgf$) \\
					\hline\hline
					\parbox{1.2cm}{Grid-\\forming} & no & \parbox[l]{3.3cm}{yes (under specific\\islanding conditions)}\\
					\hline
					\parbox{1.2cm}{Output control} & Remote setpoint control & Droop-based control\\
					\hline
		\end{tabular}}}\\
	\vspace{0pt}
		%
		\subfloat[Utility (operator) owned grid-forming DERs vs. facility level grid-forming DERs.]{\label{tab:acdc}
			\resizebox{8.7cm}{!}{
				\begin{tabular}{x{\mw}|x{\mww}|x{\mwww}}
					\hline
					Attribute & Utility-owned ($\setDERgfutil$) & Facility level  ($\setDERgffacility$) \\
					\hline\hline
					Ownership & Utility & Facility\\
					\hline
					\parbox{1.5cm}{Islanding\\condition} & \parbox[c]{4.5cm}{microgrid not connected to TN but can stay connected to other microgrid(s) (explained later in \eqref{eq:islandConditionSyn})} & \parbox[c]{4.5cm}{microgrid operates as\\an isolated island (explained later in \eqref{eq:islandConditionVsi})} \\
					\hline
		\end{tabular}}
		\vspace{0pt}
	}
\vspace{0pt}
		\caption{Comparison of DER categories.}\label{fig:derClassificationTab}
	\end{table}

	\paragraph*{DER classification}
	Our classification is based on the output behaviour and service capabilities of DERs. 
	First, we distinguish between grid-forming DERs (which provide voltage and frequency references) and DERs whose active (P) and reactive (Q) power output is controlled by PQ inverters. We denote the sets of grid-forming and PQ Inverter (PQI-) controlled DERs by $\setDERgf$ and $\setDERpq$, respectively. Then, there are a further two sub-categories of PQI-controlled DERs: those whose PQ setpoints can be remotely controlled (denoted by $\setDERpqvar$), and others whose PQ setpoints are fixed (denoted by $\setDERpqfixed$). Since the output of the DERs belonging to the set $\setDERpqfixed$ does not vary with the grid conditions, they can be considered as grid-noninteractive DERs. On the other hand, since the output of the DERs in the sets $\setDERgf$ and $\setDERpqvar$ can change with grid-conditions, we consider them as grid-interactive DERs (denoted by $\setDERgi$); see~\cref{tab:interactive}. In order to distinguish the DERs in set $\setDERpqfixed$, we refer to them as distributed generators (DGs). Since the output of these DGs cannot be changed, if operating bound violations occur, then they need to be disconnected either by remote means or through autonomous disconnections. 
	
	In contrast, the grid-interactive DERs can stay connected to the DN as zero output sources even under  fluctuations in the network state. Particularly, we assume that these DERs are fitted with low-voltage and low-frequency ride through (LVRT and LFRT) functionalities. This allows DERs to stay connected to the DN during temporary voltage and frequency bound violations at nodes. Furthermore, the output of grid-interactive DERs can be changed by two control mechanisms. In the case of grid-forming DERs ($\setDERgf$),  droop-based primary control is activated under specific islanding conditions. In the case of DERs in the set $\setDERpqvar$, their active-reactive (PQ) setpoints can be controlled by the SA system; see~\cref{tab:grid-forming}. 
	
	
	
	\def \I {\mathcal{I}}
	\def \NI {\N[\I]}
	Let $\NI \subseteq N$ denote a microgrid island within the DN. Also, let $J(S)$ be the set of DN nodes where DERs in the set $S\subseteq\setDER$ are located, with $\derNode = i$ indicating that DER $\der$ is connected to node~$i$. Recall from \cref{sec:microgridModel} that a microgrid island can consist of one or more microgrids. Based on the number of DERs contributing to  grid-forming services, a microgrid island can be in the following operating modes~\cite{microgridEmergencyControl}: 
	\begin{enumerate}
		\item Single-Master Operation (SMO): One DER operates as a single grid-forming DER (i.e. $\abs{J(\setDERgf)\bigcap\NI} = 1$), while all other DERs operate in the PQ mode. 
		\item Multi-Master Operation (MMO): More than one DER (but not necessarily all) operate as grid-forming DERs (i.e. $\abs{J(\setDERgf)\bigcap\NI} \ge 2$). 
	\end{enumerate} 
	In multi-master operation, the output of multiple grid-forming DERs changes based on nodal voltage and frequency values under the droop control constraints. These constraints ensure appropriate power sharing among DERs based on their capacities. The nodal frequencies (resp. voltages) are used for active (resp. reactive) power sharing. Our network model for radial DNs is sufficiently flexible to allow DN operations in both SMO and MMO modes.

	Finally, there are two sub-categories of grid-forming DERs, namely utility (or operator) owned grid-forming DERs and grid-forming DERs belonging to some facilities such as hospitals or other high priority loads. Each of these categories contributes to grid-forming services depending on the specific islanding conditions; see~\cref{tab:acdc}. 
	Let $\krc{\state}{\der} = 1$ if the islanding condition for DER $\der\in\setDERgf$ is satisfied, and $\krc{\state}{\der} = 0$ otherwise. 
	The two main islanding conditions of interest are as follows:
	\vspace{0.2cm}
	\begin{enumerate}[wide, labelwidth=!, labelindent=10pt]
		\item  A utility grid-forming DER contributes to grid-forming services when the node to which it belongs becomes a part of a microgrid island (i.e. the node is not connected to the TN). Consider a DER  $\der\in\setDERgfutil$ and a microgrid $\N[k]$ such that $\derNode = i \in\N[k]$. Then, DER $\der$ contributes to grid-forming  if and only if $\N[k]$ is not connected to the TN, or equivalently, at least one connecting line along the path connecting node $i$ to the substation is open, i.e. 
		\begin{equation*}
		\krc{\state}{\der} = 1 \iff \exists \ (m,n) \in \setMicrogrid \bigcap \P_i \ \text{ such that } \ \klinec{\state}{mn} = 1.
		\end{equation*}
		We formulate this condition using the following mixed-integer linear constraints:
		\begin{subequations}\label{eq:islandConditionSyn}
			\begin{alignat}{8}
			\label{eq:islandConditionSyn1}
			\krc{\state}{\der} &\ge \klinec{\state}{mn}\qquad \forall\ (m,n)\in \setMicrogrid\bigcap\P_i \\
			\label{eq:islandConditionSyn2}
			\krc{\state}{\der} &\le \ssum_{(m,n)\in(\setMicrogrid\bigcap\P_{i})}\ \ \klinec{\state}{mn}. 
			\end{alignat}
		\end{subequations}
	
		\item 
		The facility level DERs (denoted by $\setDERgffacility$) also contribute to grid-forming services when the microgrid to which they belong operates as an isolated island (i.e. not connected to the TN nor to any other microgrid). Consider a DER $\der\in\setDERgffacility$ and a microgrid $\N[i]$ such that $\derNode\in\N[i]$. Then, DER $\der$ contributes to grid-forming  if and only if all the connecting lines connecting the microgrid $\N[i]$ to the TN and other microgrids are open, i.e. 
		\begin{equation*}
		\krc{\state}{\der} = 1 \iff \klinec{\state}{mn} = 1\ \forall \ (m,n) \in \setMicrogrid_{i}.  
		\end{equation*}
		We formulate this condition using the following mixed-integer linear constraints:
		\begin{subequations}\label{eq:islandConditionVsi}
			\begin{alignat}{8}
			\label{eq:islandConditionVsi1}
			\krc{\state}{\der} &\ge \big(\ssum_{(m,n)\in\setMicrogrid_{i}}\  \klinec{\state}{mn}\big) - (\abs{\setMicrogrid_{i}} - 1) \\
			\label{eq:islandConditionVsi2}
			\krc{\state}{\der} &\le \klinec{\state}{mn}  \qquad \forall\ (m,n)\in \setMicrogrid_{i}. 
			\end{alignat}
		\end{subequations}
	\end{enumerate}
	
	\paragraph*{DER output model}
	Next, we describe the output model for the DERs. Each grid-forming DER $\der\in\setDERgf$ consists of a microsource and a storage device (batteries or flywheels)~\cite{microgridEmergencyControl}. The microsource supplies active power (quadrants I or II) in all three regimes. Thus, the output of the microsource is constrained as follows: 
	\begin{equation}
	\label{eq:nominalGeneratorConstraint}
	\Gnc{\der} \transpose{[\pnc{\state}{\der}\quad \qnc{\state}{\der}]} \le \hnc{\der}\quad \forall\ \der \in \setDERgf,
	\end{equation}
	where $\Gnc{\der}\in\R^{6\times 2}$ is a matrix and $\hnc{\der}\in \R^{6}$ is a vector that represents the polytope as shown in \Cref{fig:microsourceModel}. 
	
	For the sake of modeling simplicity, we assume that the storage device supplies active power only in the islanded regimes, whereas it consumes active power in the grid-connected regime (quadrants III and IV); see \cref{fig:resourceModel}. One justification for this restriction is that the life of a storage device significantly degrades with frequent charging/discharging cycles~\cite{microgridEmergencyControl}. Indeed, advances in storage technology make them viable sources of power supply even in the grid-connected regime. Still our modeling assumption is relevant to situations where fixed storage capacity is set aside as contingency reserve to be used in islanded regimes. 
Thus, the output of a storage device is constrained as follows: 
		\begin{alignat}{8}
\label{eq:emergenyGeneratorConstraint}
\Gec{\der} \transpose{[\pegc{\state}{\der}\quad \qegc{\state}{\der}]} + \Hec{\der} \krc{\state}{\der} &\le \hec{\der}&&\quad \forall\ \der \in \setDERgf, 
\end{alignat}
	where the $\Gec{\der},\Hec{\der}\in\R^{8\times 2}$ are matrices and $\hec{\der}\in\R^8$ is a vector such that the DER operates in  quadrants III and IV when $\krc{\state}{\der} = 0$; and in quadrants I and II  when $\krc{\state}{\der} = 1$; see~\cref{fig:storageDeviceModel}. 
	
	The total output of the DER is given by: 
	\begin{align}
	\begin{aligned}
	\prc{\state}{\der} &= \pnc{\state}{\der} + \pegc{\state}{\der} \quad&&\forall\ \der \in \setDERgf\\
	\qrc{\state}{\der} &= \qnc{\state}{\der} + \qegc{\state}{\der}&&\forall\ \der \in \setDERgf.
	\end{aligned}
	\end{align}
	On the other hand, PQI-controlled DERs ($\setDERpq$) consist only of a microsource, and do not have a storage device. Thus, their output is constrained as in \cref{fig:microsourceModel}. We can simply assume that $\forall\ \der \in \setDERpq, \  \pegc{\state}{\der} = \qegc{\state}{\der} = 0$. 
	
	
	\def \amax {3} 
	\def \bmax {3} 
	\def \lwid {5} 
	\def \lheight {4}
	\tikzstyle{nodefont} = [align=left, font=\fontsize{6pt}{9pt}\selectfont]
	\tikzstyle{qfont} = [align=center, font=\fontsize{6pt}{9pt}\selectfont]
	\tikzstyle{axis} = [thick,>=latex]

	\begin{figure}[htbp!]
		\subfloat[\footnotesize Microsource model ]{
			\label{fig:microsourceModel}
			\begin{tikzpicture}[scale=0.45]
			\draw [->,axis]  (0,0) -- (\lwid,0);
			\draw [<->,axis] (0,-\lheight) -- (0,\lheight);
			\node[nodefont, xshift=0.5cm]() at (\amax,0) {Active\\power};
			\node[nodefont, align=center, xshift=0.5cm,yshift=2]() at (0,\lheight) {Reactive power};
			
			\node[] (a1) at (0,\bmax) {}; \node[] (a2) at (\amax/2,\bmax) {}; \node[] (a3) at (\amax,\bmax/3) {}; \node[] (a4) at (\amax,-\bmax/3) {}; \node[] (a5) at (\amax/2,-\bmax) {}; \node[] (a6) at (0,-\bmax) {}; 
			\path[draw] (a1.center) -- (a2.center) -- (a3.center) -- (a4.center) -- (a5.center) -- (a6.center); 
			
			\node[nodefont,left=-0.2 of a1]() {$\qnc{max}{}$}; \node[nodefont,left=-0.2 of a6]() {$\qnc{min}{}$}; 
			\node[nodefont]() at (\amax-0.5,-0.5) {$\pnc{max}{}$};
			
			\node[qfont] at (\lwid/4,\lwid/3) {I}; \node[qfont] at (\lwid/4,-\lwid/3) {II};	
			\node[] () at (1.5*\lwid,0) {};
			\end{tikzpicture}
		}
		\subfloat[\footnotesize Storage device model ]{
			\label{fig:storageDeviceModel}
			\begin{tikzpicture}[scale=0.45]
			\draw [<->,axis]  (-\lwid,0) -- (\lwid,0);
			\draw [<->,axis] (0,-\lheight) -- (0,\lheight);
			\node[nodefont, xshift=0.5cm]() at (\amax,0) {Active\\power};
			\node[nodefont, align=center, xshift=0.5cm,yshift=2]() at (0,\lheight) {Reactive power};
			\node[qfont] at (\lwid/4,\lwid/4) {I}; \node[qfont] at (\lwid/4,-\lwid/4) {II}; \node[qfont] at (-\lwid/4,-\lwid/4) {III}; \node[qfont] at (-\lwid/4,\lwid/4) {IV};
			
			\node[] (a1) at (0,\bmax) {}; \node[] (a2) at (\amax/2,\bmax) {}; \node[] (a3) at (\amax,\bmax/3) {}; \node[] (a4) at (\amax,-\bmax/3) {}; \node[] (a5) at (\amax/2,-\bmax) {}; \node[] (a6) at (0,-\bmax) {}; \node[] (a7) at (-\amax/2,-\bmax) {}; \node[] (a8) at (-\amax,-\bmax/3) {}; \node[] (a9) at (-\amax,0) {}; \node[] (a10) at (-\amax,\bmax/3) {}; \node[] (a11) at (-\amax/2,\bmax) {}; 
			
			\path[draw] (a2.center) -- (a3.center) -- (a4.center) -- (a5.center) -- (a6.center) -- (a7.center) -- (a8.center) -- (a10.center) -- (a11.center) -- cycle; 
			\node[nodefont]() at (\amax-0.5,-0.5) {$\pegc{max}{}$}; \node[nodefont,below left=-0.2 and -0.2 of a1]() {$\qegc{max}{}$}; \node[nodefont,above left=-0.2 and -0.2 of a6]() {$\qegc{min}{}$}; \node[nodefont,below left=-0.2 and -0.2 of a9]() {$\pegc{min}{}$}; 
			\end{tikzpicture}
		}
		\caption{DER output model~\cite{shelarAminHiskens}.}\label{fig:resourceModel}
	\end{figure}
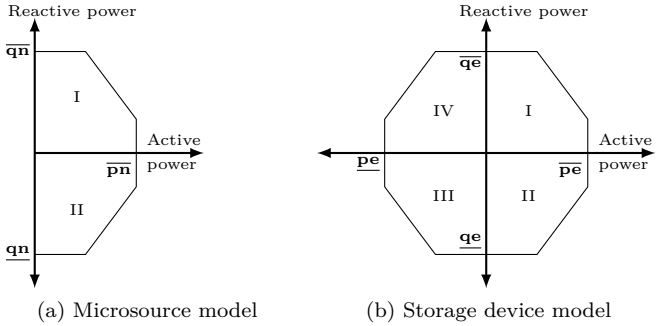
	
	\paragraph*{Droop control equations}
	We model the regulation services provided by one or more grid-forming DERs using  voltage and frequency droop control equations~\cite{islandingControlStrategies}. This allows the DERs to adjust their active and reactive power outputs based on local voltage and frequency measurements, thus  eliminating the need for explicit coordination among DERs (for the purpose of regulation). 
	
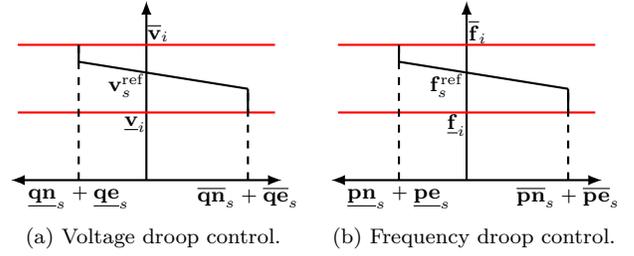
\begin{figure}[htbp!]
		\subfloat[Voltage droop control.\vspace{-8pt}]{\label{fig:droopVoltage}
			\scalefont{0.8}
			\begin{tikzpicture}[scale=0.45]
			
			\def \lwid {4}
			\def \vmax {5}; \def \vmin {3}; \def \base {1};
			\def \vone {4.5}; \def \vtwo {3.7}; 
			\draw [thick,<->,>=latex] (-\lwid,\base) -- (\lwid,\base);
			\draw [thick,->,>=latex] (0,1) -- (0,6.3);
			
			\def \lwid {3.8}
			\draw [thick,-,red] (-\lwid,\vmax) -- (\lwid,\vmax);
			\draw [thick,-,red] (-\lwid,\vmin) -- (\lwid,\vmin);
			
			\def \qrlone {-2}; 
			\def \qruone {3}; 
			
			\foreach \where/\height in {\qrlone/\vmax,\qruone/\vtwo}	
			\draw [thick,-,dashed] (\where,\base) -- (\where,\height);
			
			\draw [thick,-] (\qrlone,\vmax) -- (\qrlone,\vone);
			\draw [thick,-] (\qruone,\vtwo) -- (\qruone,\vmin);
			
			\draw [thick,-] (\qrlone,\vone) -- (\qruone,\vtwo);
			
			\foreach \where/\what in {\qrlone/$\qnc{min}{\der}+\qegc{min}{\der}$, \qruone/$\qnc{max}{\der}+\qegc{max}{\der}$}
			\node[align=center] at (\where,0.5) {\what};
			
			\node[] at (-0.33,\vmin-0.35) {$\nucmin{i}$};	
			\node[] at (-0.6,3.8) {$\nucref{\der}$};	
			\node[] at (0.35,\vmax+0.35) {$\nucmax{i}$};	
			\end{tikzpicture}
		}
		\subfloat[Frequency droop control.\vspace{-8pt}]{\label{fig:droopFrequency}
			\scalefont{0.8}
			\begin{tikzpicture}[scale=0.45]
			
			\def \lwid {4}
			\def \fmax {5}; \def \fmin {3}; \def \base {1};
			\def \fone {4.5}; \def \ftwo {3.7}; 
			\draw [thick,<->,>=latex] (-\lwid,\base) -- (\lwid,\base);
			\draw [thick,->,>=latex] (0,1) -- (0,6.3);
			
			\def \lwid {3.8}
			\draw [thick,-,red] (-\lwid,\fmax) -- (\lwid,\fmax);
			\draw [thick,-,red] (-\lwid,\fmin) -- (\lwid,\fmin);
			
			\def \prlone {-2}; 
			\def \pruone {3}; 
			
			\foreach \where/\height in {\prlone/\fmax,\pruone/\ftwo}	
			\draw [thick,-,dashed] (\where,\base) -- (\where,\height);
			
			\draw [thick,-] (\prlone,\fmax) -- (\prlone,\fone);
			\draw [thick,-] (\pruone,\ftwo) -- (\pruone,\fmin);
			
			\draw [thick,-] (\prlone,\fone) -- (\pruone,\ftwo);
			
			\foreach \where/\what in {\prlone/$\pnc{min}{\der}+\pegc{min}{\der}$, \pruone/$\pnc{max}{\der}+\pegc{max}{\der}$}
			\node[align=center] at (\where,0.5) {\what};
			
			\node[] at (-0.3,\fmin-0.37) {$\fc{min}{i}$};	
			\node[] at (-0.6,3.8) {$\fc{ref}{\der}$};	
			\node[] at (0.33,\fmax+0.39) {$\fc{max}{i}$};	
			\end{tikzpicture}
		}
		\caption[]{Droop control model~\cite{islandingControlStrategies}.}
		\label{fig:droopEquations}
	\end{figure}

	The output changes of a grid-forming DER $\der\in\setDERgf$ depend on  whether or not it is contributing to regulation  (i.e. $\krc{\state}{\der}$ = 1 or 0) based on the islanding conditions (see~\eqref{eq:islandConditionSyn} and \eqref{eq:islandConditionVsi}). Then, the classical voltage droop equation~\cite{islandingControlStrategies} can be refined to model the reactive power output of a grid-forming DER as follows (see~\cref{fig:droopVoltage}):
	\begin{equation}\label{eq:islandVoltDroop}
	\begin{split}
	\abs{\nuc{\state}{i} - \left(\nucref{\der} - \kqc{\der}(\qrc{\state}{\der} - \qrc{ref}{\der})\right)} \le \left(1-\krc{\state}{\der}\right) \bigM\\
	\forall\ \der\in\setDERgf, i\in\N \ \text{ and } \ i= \derNode.
	\end{split} 
	\end{equation}
	\cref{eq:islandVoltDroop} implies that when a DER provides regulation, it contributes more (resp. less) reactive power as the voltage drops (resp. rises) relative to a reference value. 

	Similarly, the classical frequency droop control equation~\cite{islandingControlStrategies} can be refined to model the active power output of a grid-forming DER as follows (see~\cref{fig:droopFrequency}):  
	\begin{equation}\label{eq:islandFrequencyDroop}
	\begin{split}
	\abs{\fc{\state}{i} - \left(\fc{ref}{\der} - \kpc{\der}\left(\prc{\state}{\der} - \prc{ref}{\der}\right)\right)} \le \left(1-\krc{\state}{\der}\right) \bigM\\
	\forall\ \der\in\setDERgf, i\in\N \ \text{ and } \ i= \derNode.
	\end{split}
	\end{equation}
 \cref{eq:islandFrequencyDroop} ensures proper power sharing in the sense that DERs can adjust their active power contributions for frequency regulation depending on their individual capacities. 	The reference setpoints ($\fc{ref}{\der},\nucref{\der},\prc{ref}{\der},\qrc{ref}{\der}$) and the droop coefficients ($\kpc{\der},\kqc{\der}$) are given constants.\footnote{The secondary control of voltage and frequency regulation could change the reference setpoints of the DERs, namely the voltage, frequency, active and reactive power setpoints~\cite{hierarchicalSecondaryControl2}. Secondary control may also include changing the droop coefficients of the DERs~\cite{secondaryControlDroopChange}. However, for the sake of simplicity, we  consider primary (but not secondary) control in this paper. } 

	As in \cite{part1}, we assume that each node has a DG (i.e. a grid-noninteractive DER) without loss of generality. Then, similar to the loads, we model the dependence of DG connectivity on the nodal voltage and frequency as follows:
	\begin{alignat}{8}
		\label{eq:voltageDisconnectDG2}		
	\kgc{\state}{i} &\ge \nugc{min}{i} - \nuc{\state}{i}, \quad &&
	\kgc{\state}{i} &&\ge \nuc{\state}{i} - \nugc{max}{i} \quad && \forall\ i\in \N,\\
	\label{eq:frequencyDisconnectDG} 
	\kgc{\state}{i} &\ge \fgc{min}{i} - \fc{\state}{i}, \quad && 
	\kgc{\state}{i} &&\ge \fc{\state}{i} - \fgc{max}{i} \quad && \forall\ i\in \N.
	\end{alignat}
	Eqs. \eqref{eq:voltageDisconnectDG2}-\eqref{eq:frequencyDisconnectDG} imply that a DG will disconnect if the corresponding nodal voltage or frequency violates safe operating bounds. 
	
	The net power consumed at a node $i$ is the power consumed by the load minus the power generated by the DGs and other grid-interactive DERs at that node, i.e. 
	\begin{subequations}\label{eq:netConsumptionIslanding}
		\begin{alignat}{8}
		\ptc{\state}{i} & = \pcc{\state}{i} - \pgc{\state}{i} - \ssum_{\der \in\setDERgi | \derNode = i} \ \ \prc{\state}{\der}  \qquad && \forall\ i\in\N\\
		\qtc{\state}{i} & = \qcc{\state}{i} - \qgc{\state}{i} - \ssum_{\der \in\setDERgi |  \derNode = i} \ \  \qrc{\state}{\der}  && \forall\ i\in\N.
		\end{alignat}
	\end{subequations}
	
	
	Finally, we summarize the LinDistFlow and connectivity constraints described in \cite{part1} as follows:
	\begin{alignat}{4}
	\label{eq:conserveRealApprox2} \Pc{\state}{ij} &= \sum_{k:(j,k)
		\in\E} \Pc{\state}{jk} + \ptc{\state}{j} \qquad\qquad && \forall\ (i,j)\in\E 	\\
	\label{eq:conserveReactiveApprox2} \Qc{\state}{ij} &= \sum_{k:(j,k)
		\in\E}\Qc{\state}{jk} + \qtc{\state}{j} && \forall\ (i,j)\in\E	\\
	\label{eq:dgConnectedConstraint2}
	\pgc{\state}{i} &= \left(1-\kgc{\state}{i}\right) \pgc{max}{i} && \forall\ i \in \N\\
	\qgc{\state}{i} &= \left(1-\kgc{\state}{i}\right) \qgc{max}{i} &&\forall\ i \in \N\\
	\label{eq:loadControlEquation2} 
	\pcc{\state}{i} &=  \lcc{\state}{i}\pcc{max}{i}, \quad
	\qcc{\state}{i} =  \lcc{\state}{i}\qcc{max}{i}  &&\forall\ i \in \N\\
	\label{eq:loadControlParameterConstraint2}
	(1&-\kcc{\state}{i})\lcc{min}{i} \le  \lcc{\state}{i} \le  \left(1-\kcc{\state}{i}\right) &&\forall \ i\in\N \\
	\label{eq:voltageDisconnectLoads2}		
	\kcc{\state}{i} &\ge \nucc{min}{i} - \nuc{\state}{i}, \quad 
	\kcc{\state}{i} \ge \nuc{\state}{i} - \nucc{max}{i}  && \forall\ i\in \N.
	\end{alignat}
	This completes the discussion of our multi-regime microgrid network model with parallel operation of DERs. 

	\section{Bilevel Optimization Problem}\label{sec:bilevel}
	In \cite{part1}, we modeled the sequential interaction between the attacker and operator as a bilevel mixed-integer problem (BiMIP). We now extend this model to include microgrid operations and DER dispatch capabilities. Our revised BiMIP formulation considers multi-regime microgrid operations with multiple DERs/DGs. It also accounts for TN-side  voltage and frequency disturbances as part of the overall disturbance model. 
	
\paragraph*{TN-side disruption} We consider TN-side disturbance in our attack model because the DN can face significant loss if the attacker targets the DN during an active TN failure event. In general, a TN-side disturbance (e.g. failure of a transmission line or bulk generator) can impact the system frequency as well as the substation voltage of the DN, and this can influence the attacker's strategy. We model the impact of a TN-side failure as a perturbation in the substation voltage and frequency, denoted $\vdc{}{0}$ and $\fdc{}{0}$, respectively. Then, the voltage and frequency at the substation node in the post-contingency stage can be written:
	\begin{alignat}{8}
	\label{eq:postContingencyVoltage} \nuc{\post}{0} &= \nuc{nom}{} - \vdc{}{0}, \\
	\label{eq:postContingencyFrequency} \fc{\post}{0} &= \fcnom - \fdc{}{0}. 
	\end{alignat}
	\paragraph*{DN-side disruption} For the sake of consistency, we consider the same model of DN-side disruption as in \cite{part1}, i.e. an attacker-induced compromise of the DG management system (DGMS) results in simultaneous disruption of multiple DGs. We model this attack as follows: 
	\begin{equation}
	\label{eq:dgConnectivityPostContingency}
	\kgc{\post}{i} \ge \second_i \quad \forall\quad i\in \N. 
	\end{equation}
	Let $\arcm$ denote the maximum number of DGs that the attacker can disrupt. Then, the set of all possible attacker strategies, denoted $\SecondMicro$, is given by 
	\begin{equation*}
		\SecondMicro = \{\second\in\{0,1\}^{\N} \ |\ \ssum_{i\in\N} \second_i \le \arcm\}. 
	\end{equation*}
	
	Unlike DGs (set $\setDERpqfixed$), the output of grid-interactive DERs (set $\setDERgi$) changes depending on the grid conditions. In particular, the DER output either changes autonomously based on the droop control equations, or the DERs are explicitly coordinated by the SA. The DERs are not vulnerable under our assumed disruption model because they are not affected by the compromised DGMS.

	Note that the above-mentioned disruption model can be extended to other types of attacks, including disruption of loads or circuit breakers. One can model such attacks as follows:  
	\begin{align*}
		\kcc{\post}{i}& \ge dc_i \quad &&\forall \ i \in\N\\
		\klinec{\post}{ij} &\ge dl_{ij} \ &&\forall \ (i,j) \in\E,
	\end{align*}
	where $dc\in\{0,1\}^{\N}$ and $dl\in\{0,1\}^{\E}$ denote the corresponding attacks for loads and DN lines, respectively. Thus, despite its simplicity, our approach to modeling DN-side disruptions can be applied to capture the physical impact of a broad class of security failure scenarios. This class includes Distributed Denial-of-Service (DDoS) attacks on the power grid components that can result in simultaneous failures~\cite{distributedAttackIOTDvorkin,distributedAttackIOTSoltan,loadDistributionAttack1}. Another relevant attack scenario is motivated by the vulnerabilities of Internet connected customer-side devices (e.g. smart inverters, air conditioners, water heaters), also known as Internet-of-Things (IoT) devices~\cite{distributedAttackIOTDvorkin}. An adversary can hack into these components via a cyberattack, create an IoT botnet, and can access them via the internet. Indeed, recent work in cyber-security of power systems has identified risk of correlated failures (e.g. simultaneous on/off events) induced/caused by IoT botnets~\cite{distributedAttackIOTSoltan}. In our disruption model, the impact of such an attack can be straightforwardly modeled by load/DG/line disconnects, leading to a  sudden supply-demand disturbance. However, a single point of failure such as a cyberattack on the DGMS is perhaps a more critical threat to DNs with significant penetration of DGs. 
	
\begin{remark}
	Another attack model that is well-studied in the literature considers false-data injection attacks to a (small) subset of sensors in order to inject biases in state estimates, while being undetected by anamoly detectors~\cite{Liu09falsedata,poolla,tabuada}. Available results include identification and security of \enquote{critical} sensors and attack-resilient state estimation. However, a less commonly studied aspect is that of incorrect control actions that could be implemented as a result of biased state estimation. Based on our previous work~\cite{shelarAminSunZonouz}, one can argue that our disruption model can be tailored to capture the changes in supply/demand of network nodes due to disruption of DGs/loads and/or component disconnect actions that may be induced by successful false-data injection attacks on sensor data used by the control center. 
\end{remark}
	
\begin{remark}
	Our disruption model can be extended to the compromise of grid-interactive DERs as well; see, for example, \cite{shelarAminTCNS} in which DERs in $\setDERpqvar$ are compromised by setpoint manipulation. 
\end{remark}	
	\paragraph*{Operator response model}
Recall the response capabilities (a), (b), (c) and (d) from \cref{sec:introduction2}.	Since our attack model is concerned with compromised DGMS, we rule out response (a) as an operator response. We considered (b) and (c) in \cite{part1}; see~\cref{fig:resilienceDefinition2}. 
	Our underlying assumption is that (c) is not prone to cyberattacks, because distribution utilities are being regulated under NERC CIP standards~\cite{nerccip}, which provide specific guidelines for secure \emph{reperimeterisation} of the substation cyber infrastructure. We  consider the response (d) to be executed by the SA, and thus assume that it is also secure. 
	
	The responses (b) and (c) do not consider grid-interactive DERs nor microgrid islanding capabilities.	In contrast, (d) utilizes both these capabilities, in addition to load control and preemptive disconnection of components. Particularly, we model the operator response (d) as follows: $\third \coloneqq \left(\klinec{}{}, \krc{}{}, \prc{}{},\qrc{}{}, \lcc{}{}, \kcc{}{}, \kgc{}{}\right)$. Then, the set of all response  strategies, denoted $\ThirdMicro$, can be defined as 
	$\ThirdMicro \coloneqq \{0,1\}^{\setMicrogrid} \times\{0,1\}^{\setDERgf} \times  (\R\times\R)^{\setDERgi}\times \setLoadControl \times \{0,1\}^{\N} \times \{0,1\}^{\N}$. Moreover, given the attacker-induced disruption $\second$, let the set $\ThirdMicro(\second) \coloneqq \{\third\in\ThirdMicro \ | \  \eqref{eq:dgConnectivityPostContingency} \text{ holds}\}$ denote the set of feasible response strategies available to the operator after the disruption. 
	
	For the sake of simplicity, we consider that in the pre-contingency stage, the DN is in the grid-connected regime and all components are connected. That is, there are no microgrid islands ($\klinec{\pre}{} = \zero$), and all the loads and DGs are connected to the DN ($\kcc{\pre}{} = \zero \text{ and } \kgc{\pre}{} = \zero$). Consequently, the grid-forming DERs are not contributing to regulation in the pre-contingency stage $\pre$, i.e. $\krc{\pre}{\der} = 0$ for all $\der\in\setDER$. We also assume the output of the grid-interactive DERs in mode $\pre$  to be zero, i.e. $\prc{\pre}{\der} = \qrc{\pre}{\der} = 0$ for all $\der\in\setDERgi$. These are not restrictive assumptions, however they allow us to straightforwardly compare the effectiveness of each of the response (b), (c) and (d).  
	
	\paragraph*{Post-contingency costs} 
	\label{subsec:losses}
	
	The post-contingency loss incurred by the operator, denoted $\costMicro$, is the sum of the following costs: (i) cost due to loss of voltage and frequency regulation, (ii) cost of load control, (iii) cost of load shedding, and (iv) cost of islanding:
	\begin{align}\label{eq:costGCregime2}
		\begin{aligned}
			\costMicro = &\Clovr\linfinityNorm{\nuc{nom}{} - \nuc{}{}} + \Clofr\linfinityNorm{\fcnom - \fc{}{}} \\
			&+  \Cload \ssum_{i\in\N}\ \left(\unity-\lcc{}{i}\right)  \pcc{max}{i}\\ 
			&+\left(\Cshed-\Cload\right) \ssum_{i\in\N}\ \kcc{}{i}\pcc{max}{i} \\
			&+\Cmicro \ssum_{(i,j)\in\M} \klinec{}{ij},
		\end{aligned} 
	\end{align}
	where $\Cload \in \R_{+}$ denotes the cost of per unit load control, $\Cshed \in \R_{+}$ and $\Cshed \ge \Cload$ is the cost of per unit load shed, $\Cmicro$ is the cost of a single islanding control action, $\Clovr \in \R_{+}$ is the cost of the largest deviation of nodal voltage from the nominal value $\nuc{nom}{} = 1$~pu, and $\Clofr$ is the cost of the largest deviation of nodal frequency from the nominal value $\fcnom = 1$~pu.

For a given operator response $\third\in\ThirdMicro$, let $\Xcmicro{}{}(\third)$ denote the set of post-contingency states $\xc{}{}$ that satisfy the constraints  \eqref{eq:islandingCapacity}-\eqref{eq:postContingencyFrequency}. 
	Then, we can restate our bilevel formulation  \eqref{eq:maxMinMicroGeneric} as: 
	\begin{align}\label{eq:Mm-Islanding}\tag{P-MG}
	\begin{aligned}
	\lossMicro\; \coloneqq\; &&& \max_{\second \in\SecondMicro}  \; \min_{\third\in\ThirdMicro(\second)} \;  \costMicro(\third,\xc{\post}{})  \\
	&&& \hspace{1.2cm} \text{s.t. } \xc{\post}{} \in \Xcmicro{}{}(\uc{}{}).
	\end{aligned}
	\end{align} 
Since \eqref{eq:Mm-Islanding} is a BiMIP with the same mathematical structure as the BiMIP in \cite{part1}, we solve it using the Benders Decomposition algorithm that we developed in \cite{part1}.

	\section{Computational study}
	Now, we present computational results to: (i) compare the output value of our BD algorithm with the optimal value (generated for small networks by simple  enumeration); (ii) compare the DN resilience under response capabilities (b), (c) and (d); and (iii) show the scalability of our approach to realistically large DN network sizes $\NN\in\{24,36,118\}$. 	
	
	\paragraph*{Setup for computational study} 
	\def \hp {\alpha}
	We consider three networks: modified IEEE 24-, 36-, and 118-node networks; see~\cref{fig:testNetworks} in the Appendix. The set of connecting lines $\setMicrogrid$ are shown with thick edges. The individual microgrid networks $\N[1],\cdots,\N[\abs{\setMicrogrid}]$ can be obtained by setting $\klinec{}{ij} = 1\ \forall\ (i,j)\in\setMicrogrid$. Each line $(i,j)\in\E[]$ has an identical impedance of $\resistance{ij} = 0.01, \reactance{ij} = 0.02$. Half of the nodes have a DG each and half  have a load each. 
	Consider a parameter  $\hp = \frac{6}{\NN}$. Before the contingency, each DG has active power output of $\pgc{max}{i} = \hp$, and each load has a demand of $\pcc{max}{i} = 1.25\hp$. The voltage bounds are $\nucc{min}{i} = 0.9$, $\nucc{max}{i} = 1.1$, $\nugc{min}{i} = 0.92$ and $\nugc{max}{i} = 1.08$. The reactive power values are chosen to be exactly one third that of the corresponding active power value, i.e. a 0.95 lagging power factor for each load and DG. The values are chosen such that the total net active power demand in the DN is 0.75 pu, and the lowest voltage in the network before any contingency is  close to $\nugc{min}{}$. The maximum load control parameter is $\lcc{min}{i} = 0.8$, i.e. at most 20\% of each load demand can be curtailed. For the sake of simplicity, we assume that all DGs and loads are homogeneous. The values of cost coefficients are chosen to be ${\Cload} = 100/\pcc{max}{i}, \Clovr = 100, \Clofr = 100, {\Cshed} = 1000/\pcc{max}{i}, \Cmicro = 400$. Each microgrid has one utility grid-forming (GF) DER and one facility-GF DER. Consider a parameter $\gamma = (\sum_{i\in\N}\pcc{max}{i})/(8\abs{\setMicrogrid})$. Then, each facility level DER has the following parameters: $\forall \ \der\in\setDERgffacility, \pnc{max}{\der} = \pegc{max}{\der} = \qnc{max}{\der} = \qegc{max}{\der} = \gamma$, $\kpc{\der} = 0.02$, $\kqc{\der} = 0.04$; and, each utility-owned DER has the following parameters: $\forall \ \der\in\setDERgfutil,\ \pnc{max}{\der} = \pegc{max}{\der} = \qnc{max}{\der} = \qegc{max}{\der} = 2\gamma$, $\kpc{\der} = 0.1$, $\kqc{\der} = 0.2$. These parameters are chosen such that the total capacity of grid-noninteractive DGs is 80\% of the total demand, whereas the total capacity of all grid-interactive DERs  is 75\% of the total demand of all loads. However, the total capacity of grid-interactive DERs may not be fully available to meet the demand because the microgrids are typically not of exact uniform size and topology, and the storage devices supply power only under the specific islanding configurations.

	\def \www {4.1cm}
	\def \hhh {\www}
	\paragraph*{Benders Decomposition vs. Simple Enumeration} We evaluate the ability of our implementation of the BD algorithm to compute optimal attacks in the islanding regime for small ($\NN \in \{24,36\}$) networks. For each possible cardinality of attack we first compute the optimal attack  with maximum loss using simple enumeration. Then we fix the maximum loss as $\ltarget$ for the BD algorithm. If the BD algorithm can find an attack with the same cardinality, then indeed the BD algorithm has computed the optimal attack. Otherwise, it has computed a suboptimal attack. 
		
		\def \hh {3.2cm}
		\ifNotArxivVersion
			\def \hh {3.5cm}
		\fi 
	
	\iftcnsVersion
		\def \hh {4.1cm}
	\fi
	\begin{figure}[htbp!]
		\subfloat[$\NN = 24$]{\includegraphics[width=\www,height=\hh]{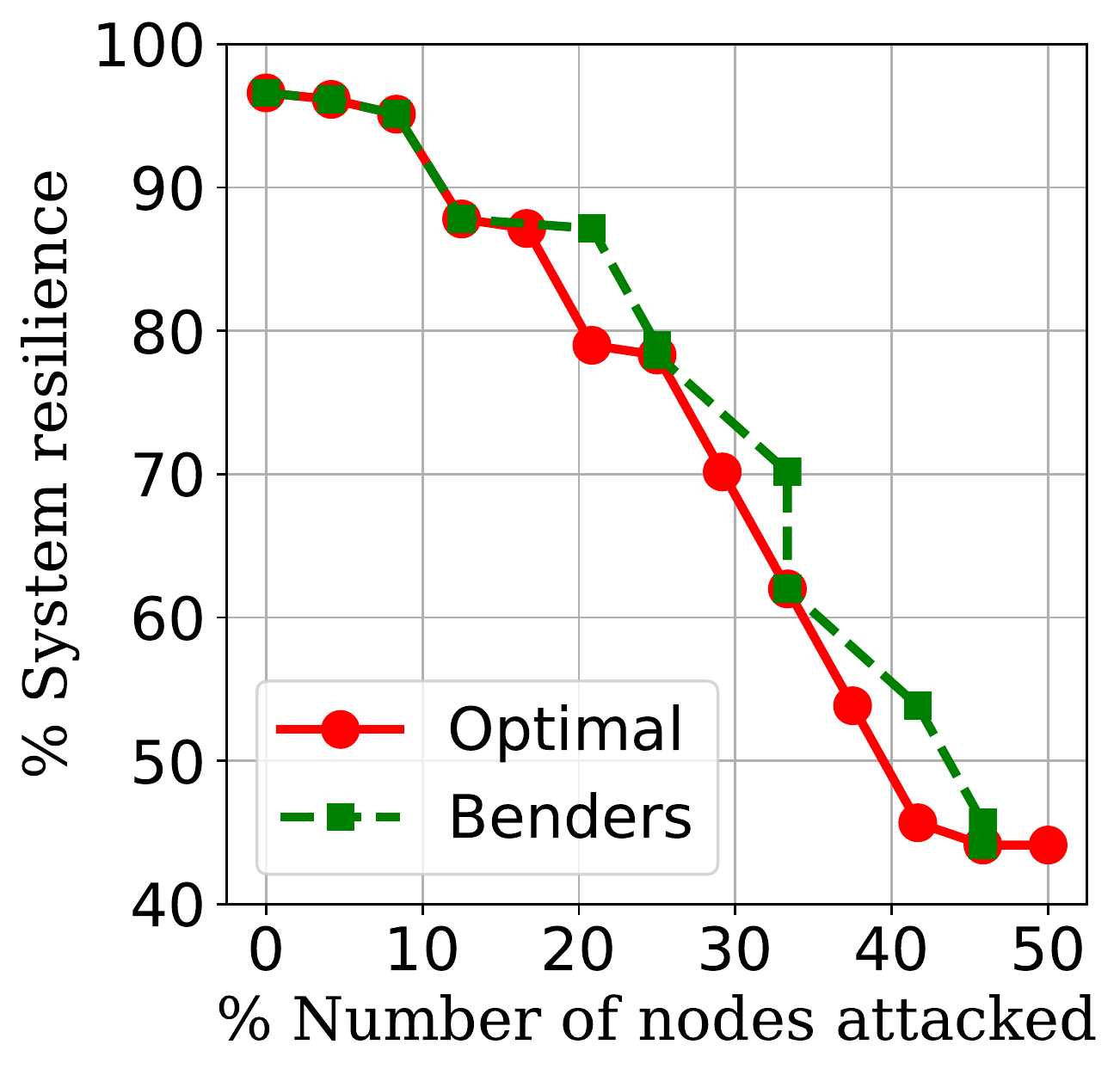}\label{fig:bruteVsBendersMicrogridN24}}
		\subfloat[$\NN = 36$]{\includegraphics[width=\www,height=\hh]{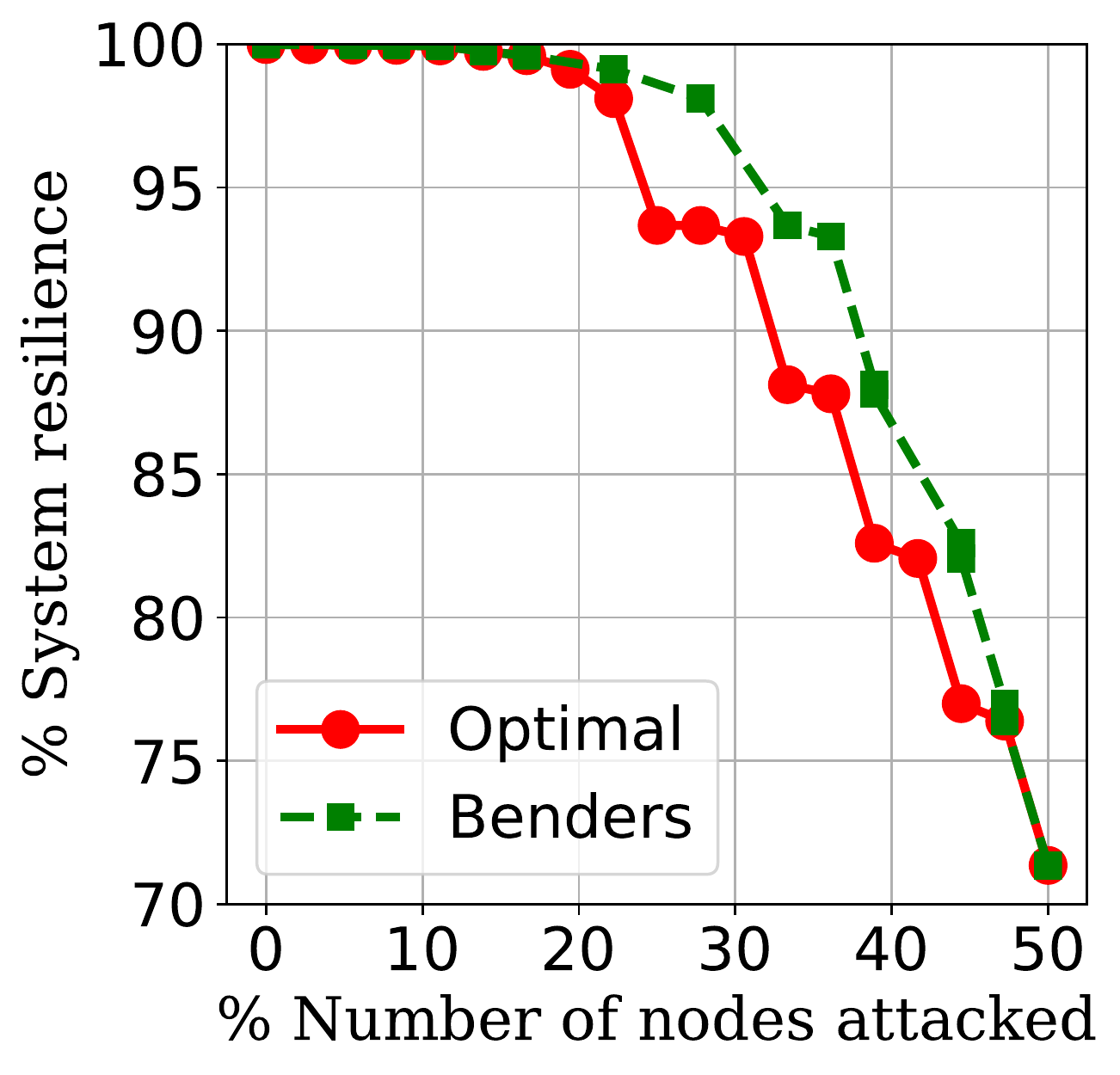}\label{fig:bruteVsBendersMicrogridN36}}
		\caption{System resilience ($\resilience = 100\left(1-\sfrac{\loss}{\lcompleteShed}\right)$) vs. $\arcm$. Near-optimal performance of BD algorithm.}
		\label{fig:bruteVsBendersMicrogrid}
	\end{figure}

	The results of the BD algorithm implemented for solving \eqref{eq:Mm-Islanding} are shown in \Cref{fig:bruteVsBendersMicrogrid}. Naturally, the attack cardinality computed by BD algorithm is greater than or equal to the optimal min-cardinality computed using simple enumeration. In some cases, however, the BD algorithm does not obtain the optimal attack. The BD algorithm involves iteratively eliminating sub-optimal attacks using Benders cuts~\cite{part1}. Each cut involved an $\epsilon$ which results in a tradeoff between the accuracy and computational time. For a very small choice of $\epsilon$, the BD algorithm eliminates exactly one sub-optimal attack in each iteration, and performs as worse as simple enumeration. For a large value of $\epsilon$, relatively more attacks,  including optimal attacks are eliminated. Hence, the BD algorithm terminates faster although with some loss of optimality. Still, for both 24- and 36-node networks, the BD algorithm computes  attacks whose cardinalities are at most 8-23\% more than the cardinalities of the corresponding optimal attacks. 
	
\begin{figure}[htbp!]
		\def \hhh {\www}
		\def \swidth {3}
		\def \drawgrid {\draw[step=1,gray,  draw opacity = 0.5] (0,0) grid (\swidth,\swidth);}
		\tikzstyle{nomargin} = [inner sep=0, outer sep = 0]
		\subfloat[$\Delta \mathrm{v}_0 = 0$\vspace{-8pt}]{
			\begin{tikzpicture}
			\node[nomargin] at (\swidth/2,\swidth/2) {\label{fig:noResponseVsSequentialVsIslandingDv0}
				\includegraphics[width=\www,height=\hhh]{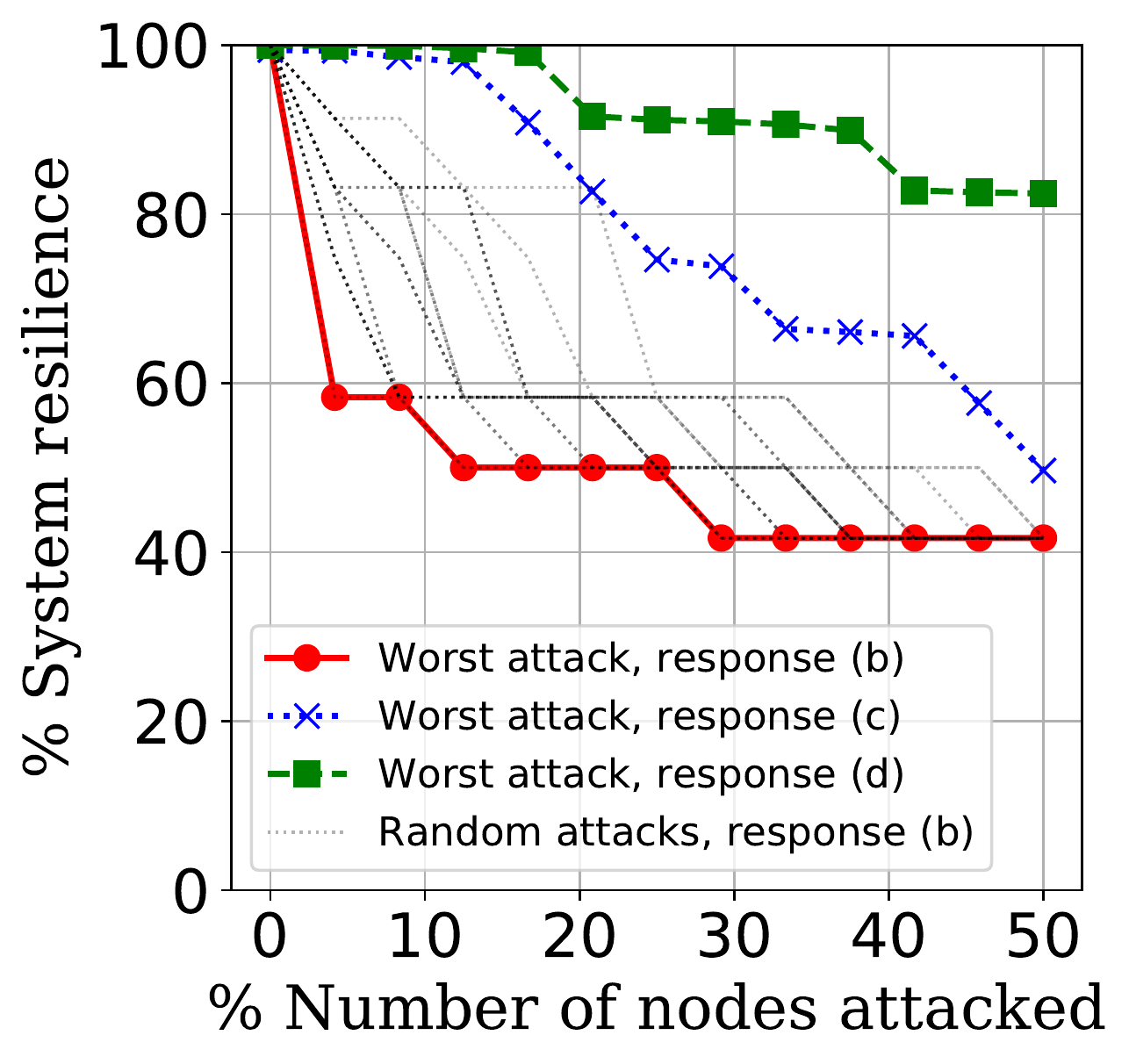}
			};
			\draw[blue!70,<->,> = latex, line width=1pt] (2.5,1.5) -- (2.5,3);
			
			\end{tikzpicture}		
		}
		\subfloat[$\Delta \mathrm{v}_0 = 0.02$\vspace{-8pt}]{
			\begin{tikzpicture}
			\node[nomargin] at (\swidth/2,\swidth/2) {\label{fig:noResponseVsSequentialVsIslandingDv2}
				\includegraphics[width=\www,height=\hhh]{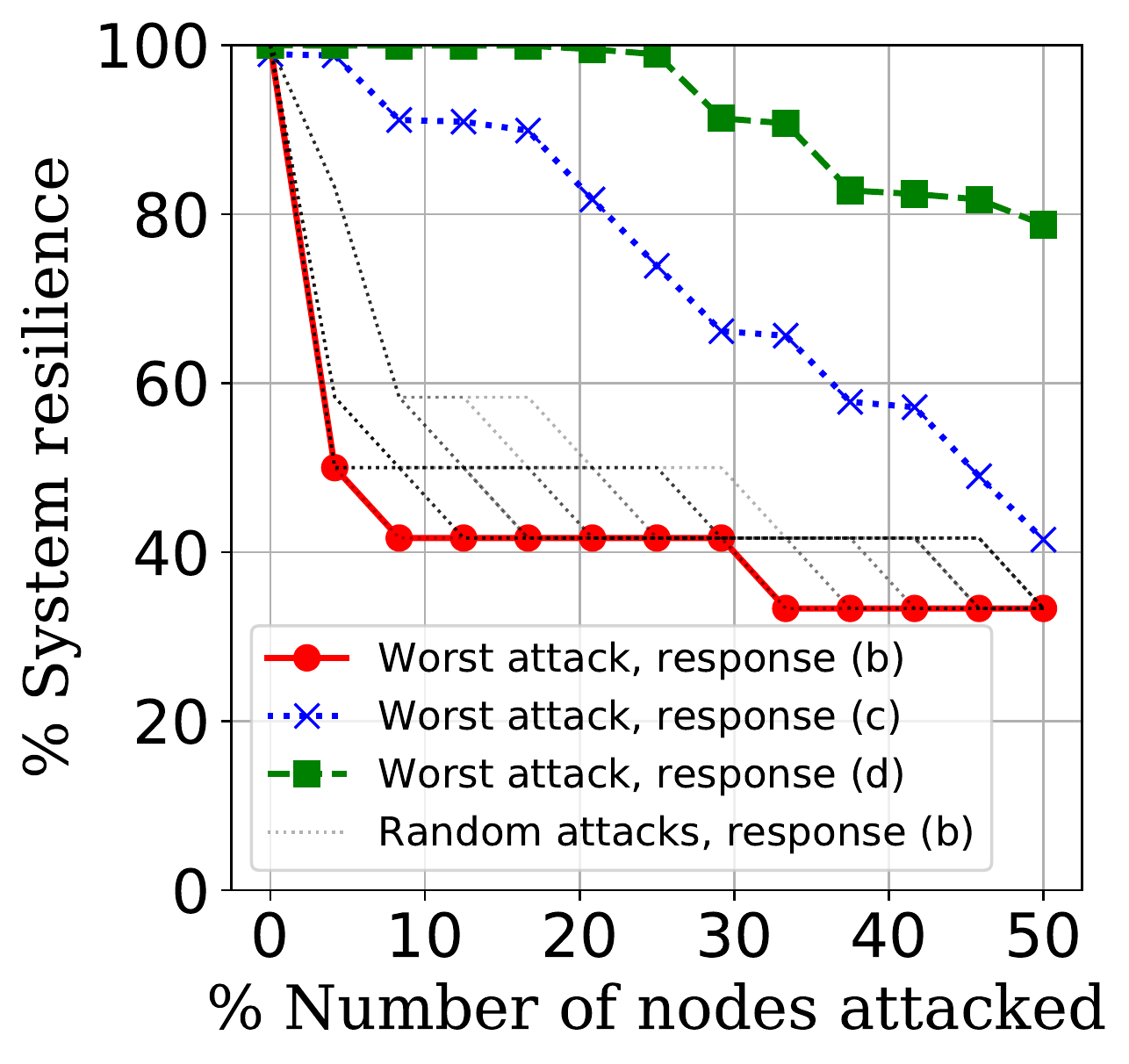}
			};
			\draw[blue!70,<->,> = latex, line width=1pt] (1.5,1.5) -- (1.5,3.3);
			
			\end{tikzpicture}		
		}
		\caption{DN resilience under varying  attacker-operator interaction scenarios. (The blue double-sided arrows indicate the value of timely microgrid response relative to the autonomous disconnections case.)}
		\label{fig:noResponseVsSequentialVsIslanding}
	\end{figure}

	\paragraph*{Value of timely response} 
	In~\cite{part1}, we used post-contingency loss to define the metric of  resilience for autonomous disconnections ($\resilienceNoResponse$) and operator response without microgrid capabilities ($\resilienceMaxmin$). In \cref{sec:introduction2}, we introduced an analogously defined metric of resilience for operator response involving microgrid islanding and DER dispatch capabilities ($\resilienceMicro$).
	\Cref{fig:noResponseVsSequentialVsIslanding} compares the resiliency values for the three cases for varying attack cardinalities, where computation of $\resilienceMicro$ and $\resilienceMaxmin$ involves using the BD algorithm to solve the corresponding BiMIPs, and $\resilienceNoResponse$ is computed using Algorithm \enquote{Uncontrolled cascade under autonomous disconnections (response (b))} in \cite{part1}. 
	Indeed, under response  (d), the SA triggers microgrid islanding and DER dispatch in a preemptive manner to reduce the impact of the attack. This leads to a smaller loss in comparison to using just load control and/or component disconnects (that is, response  (c)). Indeed, our computational results validate that $\resilienceMicro \ge \resilienceMaxmin \ge \resilienceNoResponse$. The difference between the dashed (green) and solid (red) curves in \Cref{fig:noResponseVsSequentialVsIslanding} indicate the value of response (d) relative to response (b). The difference between the dashed (green) and cross-marked (blue) curves indicate the relative value of timely response (d) over response (c). 

	\paragraph*{Scalability of the BD algorithm}
	We tabulate the performance of the BD algorithm in terms of its computational time and number of iterations to compute min-cardinality attacks for different network sizes and varying values of the resilience metric  $\resilienceTarget = 100\left(1-\ltarget/\lcompleteShed\right)$; see~\Cref{tab:bendersPerformance2}. We also note the cardinalities of attacks output by the BD algorithm as well as the corresponding DN resilience. 
	Note that the $\NN=118$ node network has $2^{118}$ possible configuration vectors. Still, with $\resilienceTarget = 80\%$, the BD algorithm computes an attack in $\approx$ 1 minute. In comparison, for the $\NN=36$ node network, the simple enumeration method took $\approx$ 6 hours. 
	\begin{table}[htbp!]
		\centering
		\caption[]{ Resiliency metric evaluated using the BD algorithm for 24-, 36- and 118-node networks. The realized resilience metric can significantly fall short of the target resilience metric ($\resilienceTarget = 100\left(1-\ltarget/\lcompleteShed\right)$); for example, when the attack cardinality changes from 6 to 7, the percentage resilience for the 24-node network decreases sharply from 98.91\% to 91.33\%. This means that the 24-node DN is at least 90\% (actual value 91.33\%) resilient to $\arcm=7$ cardinality attacks. }\label{tab:bendersPerformance2}
		\resizebox{1\textwidth}{!}{
			\begin{tabular}{|c|l|l|l|}
				\hline 
				\multicolumn{4}{|c|}{\parbox{1.2\textwidth}{\textbf{Entries are resilience metric of DN (in percentage), number of iterations (written in brackets), time (in seconds), attack cardinality.}}}\\
				\hline 
				$\resilienceTarget$ 
				& $\NN = 24$ & $\NN = 36$ & $\NN = 118$  \\
				\hline 
				$99$ &98.91, (15), 0.41, 6  & 98.95, (10), 0.37, 5 &  98.95, (8), 2.48, 4 \\
				\hline
				$95$  & 91.33, (16), 0.46, 7  & 94.12, (12), 0.51, 7 &  94.28, (15), 3.91, 11 \\
				\hline
				$90$ & 82.8, (18), 0.57, 9  & 88.23, (17), 0.91, 11 & 89.73, (20), 10.62, 16 \\
				\hline
				$85$  & 82.8, (18), 0.57, 9  & 81.9, (20), 1.23, 14  &  83.49, (29), 28.79, 25 \\
				\hline $80$ & 78.73, (21), 0.74, 12 & 71.46, (21), 1.75, 15 &  79.9, (40), 67.38, 36 \\
				\hline
			\end{tabular}
		}
	\end{table}

	\section{Multi-period DN restoration}\label{sec:restoration}
	We recall that the resilience of a system is related to its ability to not only minimize the impact of a disturbance, but also quickly recover from it; see~\cref{fig:resilienceDefinition2}. Our attack model assumes that a compromise of the DGMS leads to remote disconnection of multiple DGs. However, the actual functionality of disconnected DGs is not compromised. In response  (d), we consider that the SA has the ability to detect and obtain knowledge of the complete attack. Moreover, the SA can also control DG connectivity. We now discuss how the SA can restore the disrupted DGs, and bring the DN back to its nominal performance. 
	In this section, we present a simple MIP that models the process of restoring system performance. 
	

	Our model of the DN restoration process entails progressively reconnecting the disrupted DGs, and eventually restoration to the grid-connected mode of DN operation. We consider a multi-period horizon $\setPeriods=\{0,1,\cdots,\nperiod\}$ where the end of the $0^{th}$ period coincides with the time when the restoration actions begin; see~\cref{fig:resilienceDefinition2}. Let a period be denoted by $\period\in\setPeriods$, where each period $\period$ is of fixed time duration (say, a few minutes). Furthermore, the operator response at period $\period$ is denoted by $\uc{\period}{}$. Let $\nperiod = \period_{\text{res}} + 1$, where $\period_{\text{res}}\in\mathbb{Z}_+$ denotes the earliest time period when all disrupted DGs can be restored.
	
	A TN-disturbance may clear any time, before or after the DG reconnections. However, since our analysis is focused on determining worst-case resilience of the DN, we assume that the TN-side disturbance clears after the disrupted DGs are fully reconnected, i.e. beyond time $\nperiod$. Under the assumed detection and response capabilities of the SA, the restoration action begins only after $\period = 0$. Thus,  the operator control actions in the $0^{th}$ time period are the same as the initial contingency response, i.e. ${\uc{0}{} = \uc{\post}{}}$. Similarly, the operator control actions remain unchanged after the restoration actions are complete until just prior to the TN-disturbance clearing. Therefore, without loss of generality we can assume that $t_2 = \nperiod$, and restrict our analysis to the horizon $\setPeriods$. 
	
	We consider two types of constraints to model the restoration actions of the operator across time periods: monotonicity constraints and resource constraints. Consider a period $\period \in \{1,2,\cdots,\nperiod\}$. The monotonicity constraints for period $\period$ are: 
	\begin{subequations}\label{eq:continuityConditions}
		\begin{alignat}{8}
		\label{eq:lineRepairChange}
		\klinec{\period}{\edge} &\le \klinec{\period-1}{\edge} \quad\forall\ (i,j)\in\setMicrogrid,\\
		\label{eq:dgrestoration}
		\kgc{\period}{i} &\le \kgc{\period-1}{i}\quad\forall\ i\in \N. 
		\end{alignat}
	\end{subequations}
	\cref{eq:lineRepairChange} implies that during the restoration process, once a connecting line is closed, it remains closed until the  restoration process is completed. Similarly, \eqref{eq:dgrestoration} implies that a disconnected DG becomes operational after being reconnected, and then remains operational until the restoration is complete. The monotonicity constraints can be justified based on the practical consideration that frequently changing the status of connecting lines can create large fluctuations in nodal voltages and frequencies of the microgrids due to the low inertia of DERs. Moreover, the battery life of storage devices would reduce due to frequent changes from charging modes (quadrants III and IV) to discharging modes (quadrants I and II), and vice versa; see~\cref{fig:resourceModel}. 

	
	The resource constraint merely limits the number of DG reconnections. Specifically, we consider that during period $\period$, at most $\ngc{\period}{}$ DGs can be reconnected, where $\ngc{\period}{}$ denotes the restoration budget for that period:
	\begin{equation}\label{eq:maxrestoreDG}
	\ssum_{i\in\setDERpqfixed} \kgc{\period}{i} \ge \ssum_{i\in\setDERpqfixed}\kgc{\period-1}{i} - \ngc{\period}{}.
	\end{equation}
	Restrictions on the number of connecting line closing operations can be similarly considered. \cref{eq:maxrestoreDG} can also be justified in a way similar to that of the monotonicity constraints. Naturally, the operator wants to avoid a large number of simultaneous DG reconnections as that could lead to large voltage and frequency fluctuations. 
	
	We must choose $\nperiod$ large enough so that all disrupted DGs can be reconnected over the horizon $\setPeriods$. This implies $\nperiod = \min\{\period' | \sum_{\period=1}^{\period'}\ngc{\period}{} \ge \arcm\} + 1$. Also, we assumed that the TN-side disturbance ceases to exist at the last time period. We model this as:
	\begin{subequations}\label{eq:restorationTNconditions}
		\begin{alignat}{8}\small
		\nuc{\period}{0} &= \begin{cases}
		\nuc{nom}{} - \vdc{}{} \quad& \text{if } \period \ne \nperiod\\
		\nuc{nom}{} & \text{if } \period = \nperiod
		\end{cases} \\
		\fc{\period}{0} &= \begin{cases}
		\fcnom - \fdc{}{} \quad& \text{if } \period \ne \nperiod\\
		\fcnom & \text{if } \period = \nperiod.	
		\end{cases}
		\end{alignat}
	\end{subequations}
	Let $\Ycmicro{\period}{}(\uc{\period-1}{})$ denote the feasible set of response strategies for $\uc{\period}{}$, i.e. $\Ycmicro{\period}{}(\uc{\period-1}{}) = \{\uc{\period}{} \in \ThirdMicro\ | \  \text{such that } \eqref{eq:continuityConditions}-\eqref{eq:maxrestoreDG} \text{ hold}\}$. Also, given an operator response $\third\in\ThirdMicro$, let $\Xcmicro{\period}{}(\third)$ denote the set of network states $\xc{\period}{}$ which satisfy the constraints  \eqref{eq:islandingCapacity}-\eqref{eq:voltageDisconnectDG2} and \eqref{eq:restorationTNconditions}. Hence, the restoration problem can be posed as:
	\begin{align}\label{eq:restoration}\tag{P3}\small
	\begin{aligned}
	\hspace{-0.2cm}\lossres(\second) \coloneqq\;  \min_{ \{\third^\period\}_{\period\in\setPeriods}} &&& \sum_{\period\in\setPeriods} \costMicro\left(\third^\period,\xc{\period}{}\right) \\
	\text{s.t.} &&& \uc{0}{} \in \ThirdMicro(\second) \\
	&&& \uc{\period}{} \in \Ycmicro{\period}{}(\uc{\period-1}{}) && \forall\ \period = 1,\cdots,\nperiod\\
	&&& \xc{\period}{} \in \Xcmicro{\period}{}(\uc{\period}{})  && \forall\ \period = 0,\cdots,\nperiod.\\
	\end{aligned}
	\end{align}
	Problem \eqref{eq:restoration} is a Mixed-Integer Problem (MIP), and can be solved using off-the-shelf MIP solvers. However, due to the large number of binary variables, it can become computational expensive to solve for larger networks. In fact, we solve \eqref{eq:restoration} using a simple greedy algorithm; see~\cref{algo:greedyAlgorithm}. In each period, the operator simply chooses that response which minimizes the post-contingency loss during that time period subject to the monotonicity and resource constraints. \cref{algo:greedyAlgorithm} is based on the feature that the network state in any period depends only on the operator actions in that period, and the network state in the previous period.  The algorithm returns with the operator actions, resulting network state, and corresponding post-contingency loss for each time period. 
	\begin{algorithm}[htbp!]
		\caption{Greedy Algorithm}
		\begin{algorithmic}[1]\small
			\State $\uc{0}{}, \xc{0}{} \gets \argmin\limits_{\third\in\ThirdMicro(\second)} \quad  \costMicro(\third,\xc{}{}) \quad \text{s.t.} \quad \xc{}{} \in \Xc{}{}\left(\third\right)$.
			\For {$\period = 1,\cdots,\nperiod$} 
			\State $\uc{\period}{}, \xc{\period}{} \gets \argmin\limits_{\third\in \mathcal{Y}^\period_{\text{m}}(\uc{\period-1}{})} \quad  \costMicro(\third,\xc{}{}) \quad \text{s.t.} \quad \xc{}{} \in \Xcmicro{\period}{}(\third)$. 
			\State $\cost^\period \gets \costMicro(\uc{\period}{},\xc{\period}{})$
			\EndFor
			\State \Return $\{\uc{\period}{}, \xc{\period}{}, \cost^\period\}_{\period\in\setPeriods}$ %
		\end{algorithmic}
		\label{algo:greedyAlgorithm}
	\end{algorithm}
	
\begin{figure}[htbp!]
			\ifArxivVersion
				\def \hh {4.3cm} 
			\fi 	
		\subfloat[$\Delta \mathrm{v}_0 = 0.1$]{\label{fig:restoration1}	\includegraphics[width=4.3cm,height=\hh]{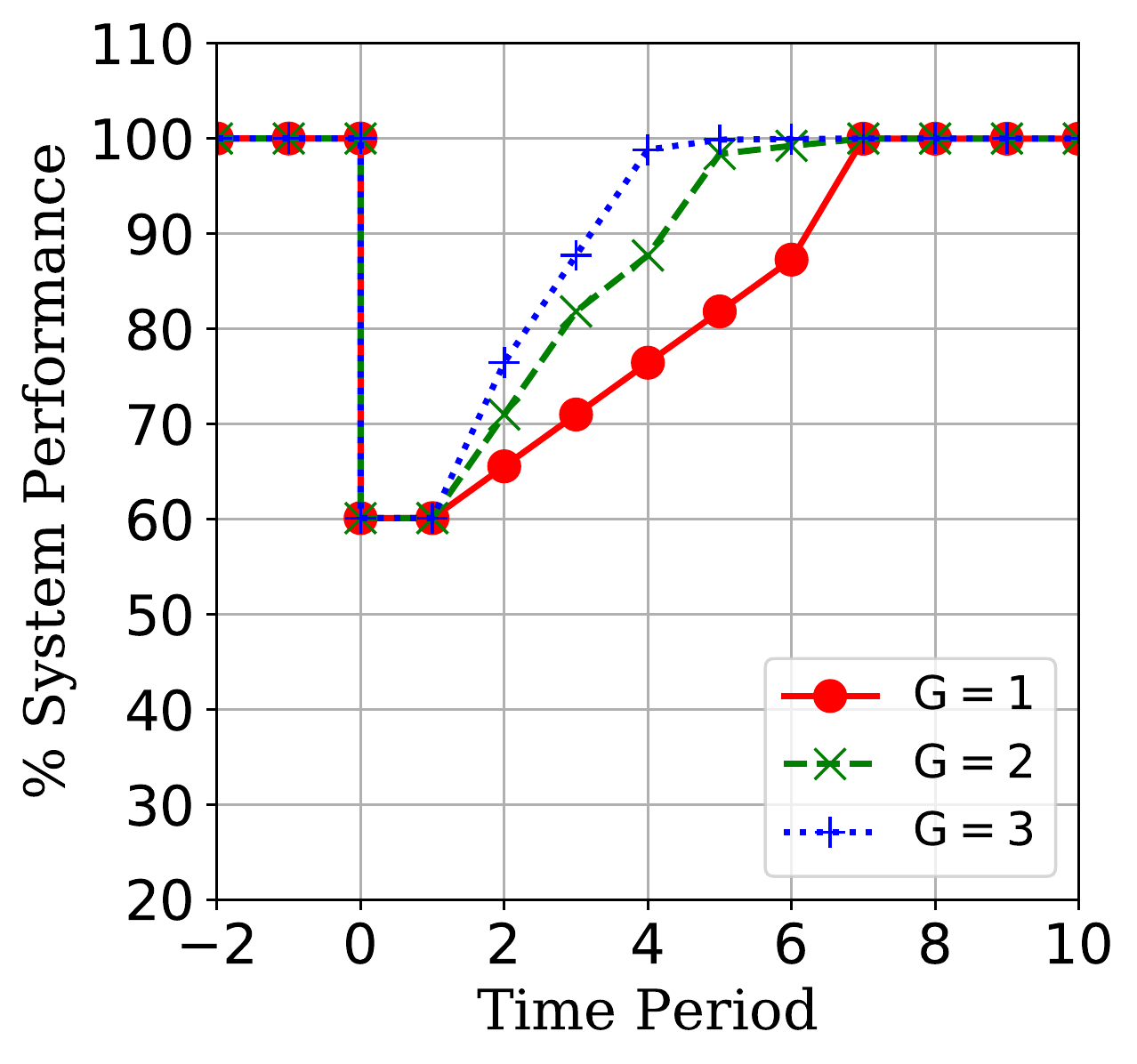}}
		\subfloat[$\Delta \mathrm{v}_0 = 0.2$]{\label{fig:restoration2}	\includegraphics[width=4.3cm,height=\hh]{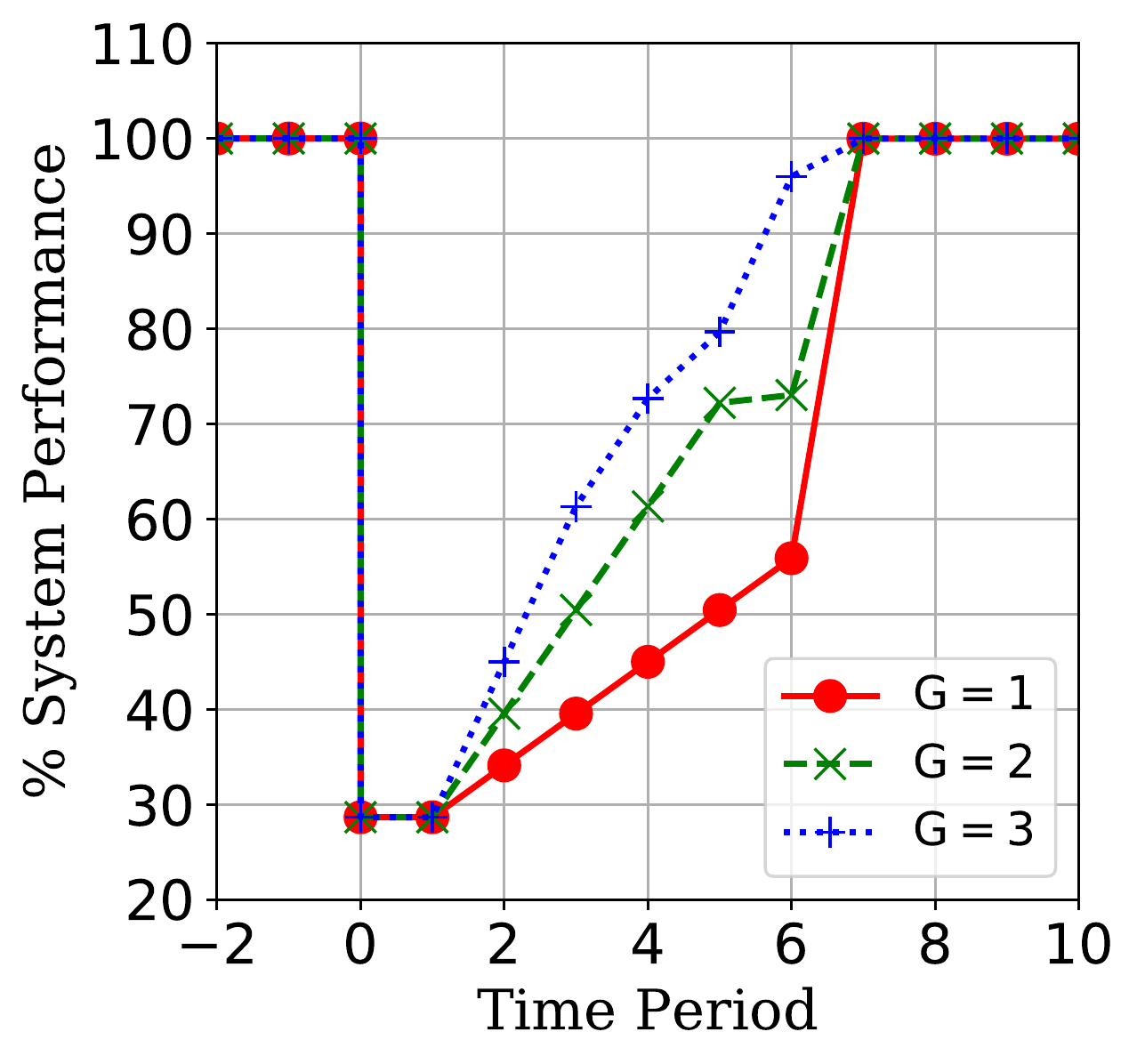}}
		\caption{Multi-period DN restoration ($\NN = 36$).}
		\label{fig:restoration}
	\end{figure}

	\Cref{fig:restoration} shows the system performance during the restoration of the DN over multiple time periods for different resource constraints. For each system restoration curve, we chose $\ngc{\period}{}$ to be a constant for all time periods $\period\in\setPeriods$. One can see that after the TN-side and DN-side  disturbances, the system performance drops. Then, as disrupted components are  connected, the system performance progressively recovers. Also, the post-contingency losses are higher for larger TN-side disturbances. However, as the restoration budget increases, the system recovers faster. 

\begin{figure}[htbp!]
		\iftcnsVersion
		\def \hh {4.3cm}
		\fi 	
		\subfloat[$\NN = 24$ \vspace{-4pt}]{\label{fig:restorationComparison1}	\includegraphics[width=4.3cm,height=\hh]{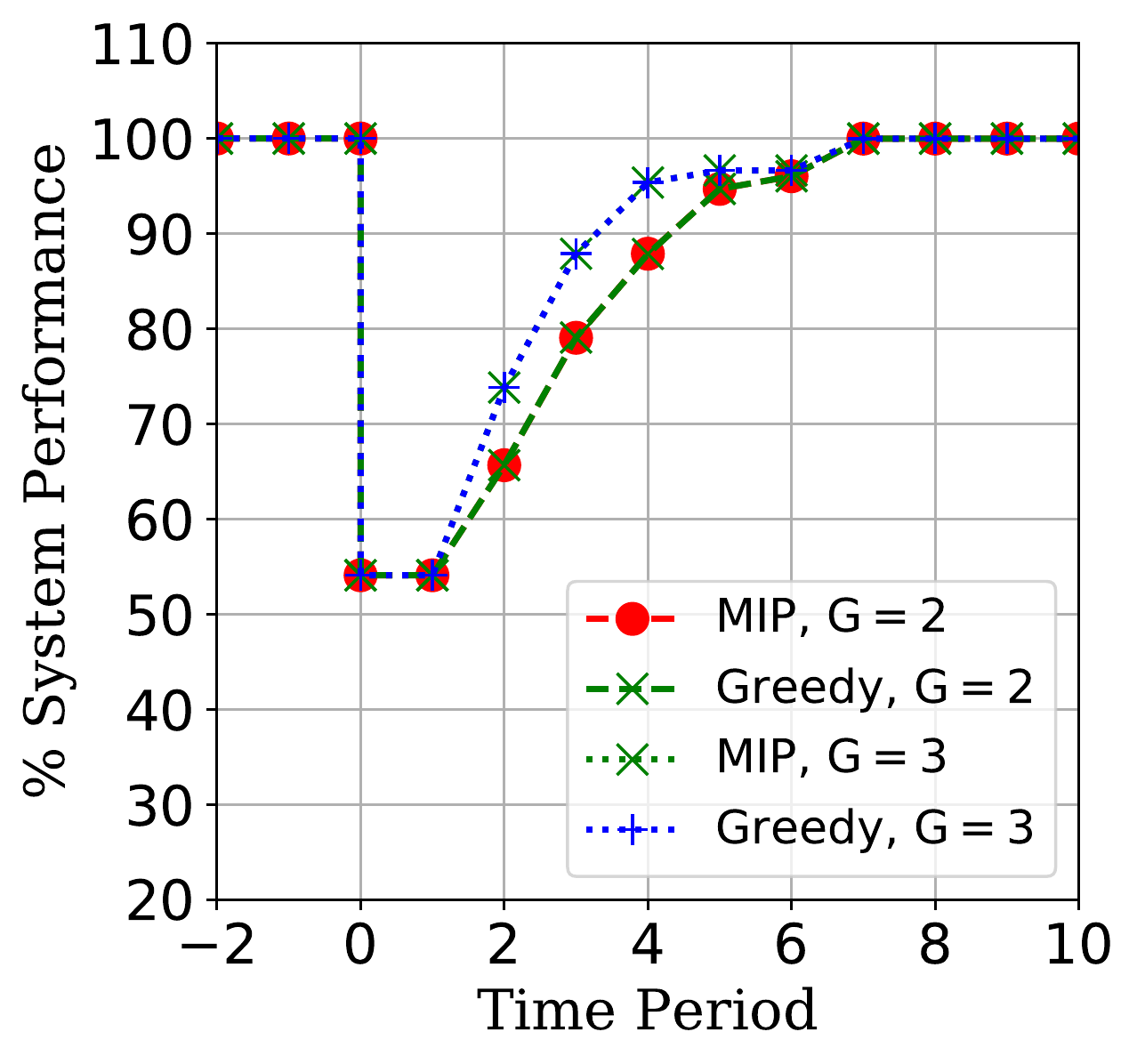}}
		\subfloat[$\NN = 36$ \vspace{-4pt}]{\label{fig:restorationComparison2}	\includegraphics[width=4.3cm,height=\hh]{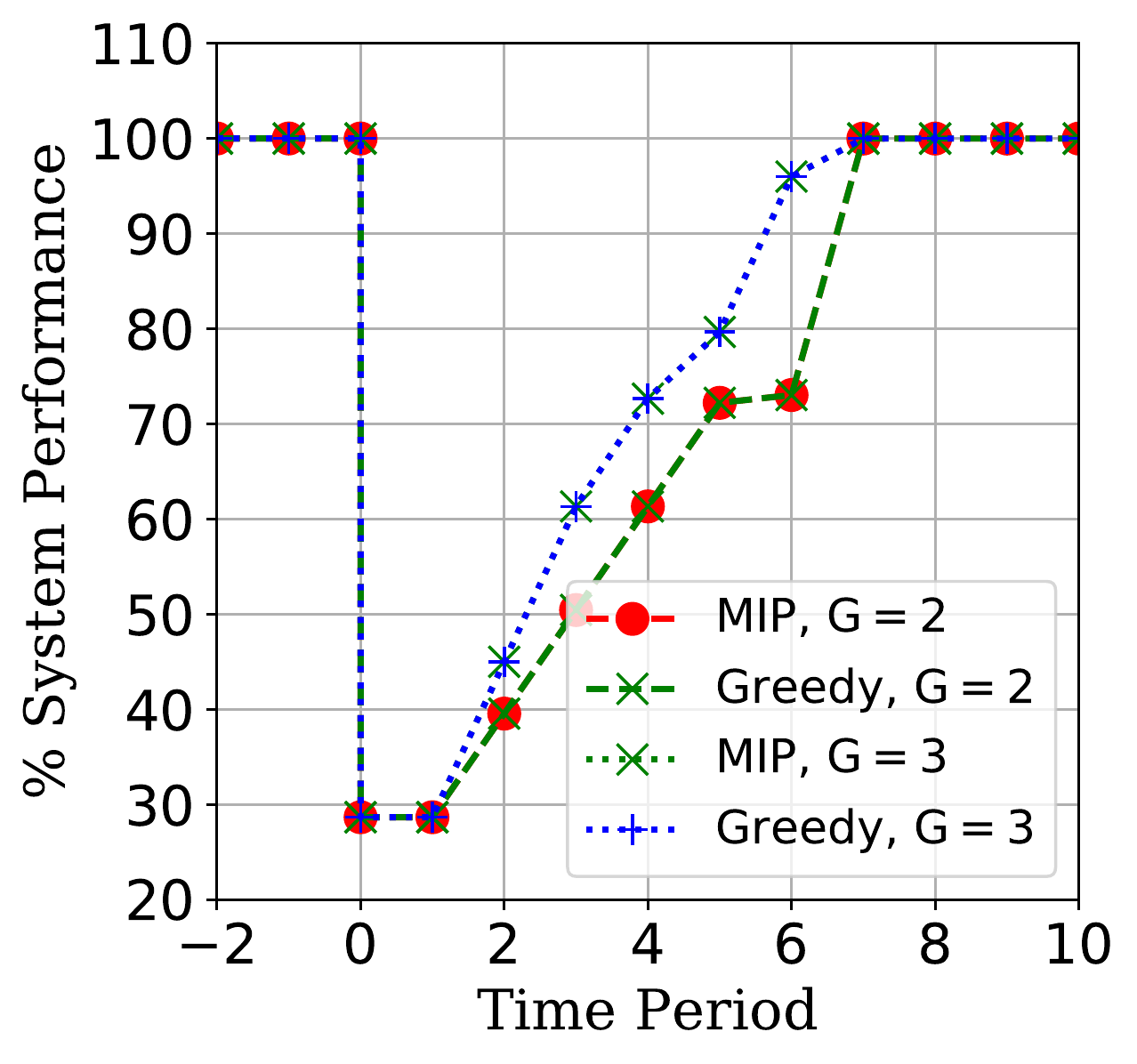}}
		\caption{Near-optimal performance of Greedy \Cref{algo:greedyAlgorithm}.} 
		\label{fig:restorationComparison}
	\end{figure}

\paragraph*{MIP vs. Greedy Recovery Algorithm}
\Cref{fig:restorationComparison} shows the comparison of the system performance recovery curves obtained using \cref{algo:greedyAlgorithm} and by directly solving the large-scale MIP for $\NN=24$ and $\NN=36$ node networks. The TN-side voltage disturbance for both the networks is set to $\vdc{}{} = 0.2$. In this experiment, we set the time limit of the (Gurobi) solver to 7200 seconds. 
While solving the large-scale MIP for  $\NN=36$ and $\ngc{}{}=3$  we were able to achieve an optimality gap of 16.54\% after 2 hours. However,  \cref{algo:greedyAlgorithm} was able to attain the same system performance recovery curve using the default solver settings (no presolve and Simplex method), and compute the near-optimal solution in approximately 10 seconds. 
	
In order to implement the response computed in \eqref{eq:restoration}, the SA may need to coordinate with the individual microgrid controllers. A detailed description of such a  communication architecture is beyond the scope of this paper. We refer the reader to \cite{hierarchicalControl} for a hierarchical control architecture which can support the coordination between SA and individual microgrid controllers. 
	
	\section{Concluding remarks}
	\label{sec:conclusions2}
	
	In this paper and its companion paper \cite{part1}, we developed a quantitative
	framework to evaluate DN resilience, which is its ability to minimize
	the impact of a disturbance, and to restore the DN back to full load
	supply. It has been shown that the impact of a broad class of
	cyberphysical failure scenarios can be modeled as DN-side disruption
	of multiple components and/or disturbances in substation voltage and
	frequency. We developed a novel network model which captures operation
	of microgrid(s) under various regimes, and single-/multi-master
	operation of DERs. Various operator response strategies were
	considered: from load control and component disconnects to microgrid
	islanding and DER dispatch. Furthermore, we formulated the
	attacker-operator interactions as bilevel mixed-integer problems, and
	developed a computational approach to efficiently solve these problems
	using a Benders decomposition algorithm. Restoration of DN performance
	over multiple time periods was considered, and a greedy algorithm was
	presented for solving this problem. Our computational results show the
	value of timely response under varying operator capabilities in
	minimizing the impact of disruption.
	
	Even though a linear power flow was used and only basic aspects of
	microgrid operations have been considered, the paper provides a rich
	and flexible modeling framework for analyzing DN resilience for more
	sophisticated attack and response capabilities. Other cyber-physical
	security scenarios can be similarly analyzed by considering a clear
	demarcation between the vulnerable and the securely controllable DN
	components. The computational approach for solving the bilevel
	formulation arising from the linear power flow approximation can, in
	principle, be extended to convex (second-order cone) relaxations of
	the nonlinear power flow model. This approach may be useful for
	optimal resource allocation~\cite{shelarAminHiskens} and security investments for DNs~\cite{shelarAminTCNS}. The framework also provides the basis for resiliency assessment of other smart infrastructure networks.

\bibliographystyle{IEEEtran}	
\iftcnsVersion
	\small
	\bibliography{IEEEfull.bib}
\fi 
\normalsize
	
\iftcnsVersion
\appendix
	\begin{figure}
		\tikzstyle{ndc} = [draw,circle,inner sep=0,minimum width=0.55cm]
		\tikzstyle{ndcSyn} = [ndc,fill,pattern=north west lines]
		\tikzstyle{ndcVsi} = [ndc,fill,pattern=vertical lines]
		\tikzstyle{mgl} = [draw,line width=3]
		\subfloat[24 node network.]{
			\resizebox{8.7cm}{!}{
				\begin{tikzpicture}
				\def \nodePar {0/1,1/2,2/3,3/4,4/5,5/6,6/7,7/8,8/9,9/10,10/11,11/12,12/13,13/14,14/15,15/16,16/17,17/18,18/19,19/20,20/21,21/22,22/23,23/24,24/25,25/26,26/27,27/28,28/29,29/30,30/31,31/32,32/33,33/34,34/35,35/36}
				
				\def \ny {0}		
				\node[ndc] (0) at (0,\ny) {0};		
				\foreach \x in {1,...,12}
				{
					\node[ndc](\x) at (\x,\ny) {\x};
				}
				
				\def \ny {1}		
				\foreach \x in {13,...,16}
				{
					\pgfmathtruncatemacro{\nx}{\x-10}
					\node[ndc](\x) at (\nx,\ny) {\x};
				}
				\def \ny {-1}		
				\foreach \x in {17,...,24}
				{
					\pgfmathtruncatemacro{\nx}{\x-13}
					\node[ndc](\x) at (\nx,\ny) {\x};
				}
				
				\def \synNodes {{1,7,13,17,20}}
				\foreach \x in {1,7,13,17,20}
				\node[ndcSyn] at  (\x) {};
				\foreach \x in {3,9,15,19,22}
				\node[ndcVsi] at  (\x) {};
				\foreach \x/\y in {0/1,1/2,2/3,3/4,4/5,5/6,6/7,7/8,8/9,9/10,10/11,11/12,2/13,13/14,14/15,15/16,3/17,17/18,18/19,6/20,20/21,21/22,22/23,23/24}
				\draw (\x) -- (\y);
				\foreach \x/\y in {0/1,6/7,2/13,3/17,6/20}
				\draw[mgl] (\x) -- (\y);
				\end{tikzpicture}}}\qquad\\
		\subfloat[36 node network.]{\resizebox{6cm}{!}{
				\begin{tikzpicture}
				\foreach \z/\x/\y in {0/0/0,1/1/0,2/2/0,16/3/0,21/4/0,22/5/0,23/6/0,24/7/0,6/2/1,7/2/2,8/3/2,9/4/2,14/1/2,15/0/3,13/1/3,10/2/3,11/3/3,12/4/3,4/1/-1,3/2/-1,5/2/-2,17/3/-1,18/3/-2,19/4/-2,20/3/-3,25/5/-1,26/5/-2,27/6/-1,28/7/-1,29/7/-2,31/6/-2,30/8/-2,32/8/-1,33/8/0,34/8/1,35/8/2,36/7/1}
				\node[ndc] (\z) at (\x,0.8*\y) {\z};
				\foreach \x in {1,6,10,16,21,27,32}
				\node[ndcSyn] at  (\x) {};
				\foreach \x in {3,8,12,18,23,29,34}
				\node[ndcVsi] at  (\x) {};
				\foreach \x/\y in {0/1,1/2,2/3,3/4,3/5,2/6,6/7,7/8,8/9,7/10,10/11,11/12,10/13,13/14,13/15,2/16,16/17,17/18,18/19,18/20,16/21,21/22,22/23,23/24,22/25,25/26,25/27,27/28,28/29,29/30,29/31,28/32,32/33,33/34,34/35,34/36}
				\draw (\x) -- (\y);
				\foreach \x/\y in {0/1,2/6,7/10,2/16,16/21,25/27,28/32}
				\draw[mgl] (\x) -- (\y);
				\end{tikzpicture}	
		}}\\
		\subfloat[118 node network.]{
			\resizebox{8.7cm}{!}{
				\begin{tikzpicture}
				
				\def \lateral(#1,#2,#3,#4) {
					\foreach \x in {#1,...,#2}
					{
						\pgfmathtruncatemacro{\nx}{\x-#3}
						\node[ndc](\x) at (\nx,0.8*#4) {\x};
					}
					\pgfmathtruncatemacro{\ns}{#1+1}
					\foreach \x in {\ns,...,#2}
					{
						\pgfmathtruncatemacro{\ef}{\x-1}
						\draw (\x) -- (\ef);				
					}
				}
				
				
				\lateral(18,27,16,7)
				\lateral(10,17,9,6)
				\lateral(4,9,2,5)
				\node[ndc] (2) at (1,5*0.8) {2};
				\node[ndc] (3) at (0,5*0.8) {3};
				\lateral(38,46,35,4)
				\lateral(28,35,26,3)
				\lateral(47,54,42,2)
				\lateral(36,37,33,2)
				\lateral(55,62,53,1)
				\lateral(96,99,92,0)
				\lateral(0,1,0,0)
				\lateral(89,95,87,-1)
				\lateral(63,69,62,-2)
				\lateral(70,77,66,-3)
				\lateral(78,80,77,-3)
				\lateral(81,85,77,-4)
				\lateral(86,88,85,-4)
				\lateral(100,106,100,-5)
				\lateral(107,113,102,-6)
				\lateral(114,118,114,-6)
				
				\foreach \x in {1,10,18,30,38,47,55,63,70,80,89,100,101,107}
				\node[ndcSyn] at  (\x) {};
				\foreach \x in {3,12,20,32,40,49,57,65,72,82,91,115,103,109}
				\node[ndcVsi] at  (\x) {};
				\foreach \x/\y in {1/2,2/3,2/4,4/28,30/36,80/81,100/114,64/78,79/86,65/89,91/96}
				\draw (\x) -- (\y);
				\foreach \x/\y in {0/1,2/10,11/18,29/30,29/38,35/47,29/55,1/63,65/89,69/70,79/80,106/107,1/100,100/101}
				\draw[mgl] (\x) -- (\y);
				\end{tikzpicture}}}
		\caption{Modified IEEE test networks. Connecting lines are indicated by thick edges. Utility-owned (resp. facility level) grid-forming DERs are indicated by northwest (resp. vertical) lines.}\label{fig:testNetworks}
	\end{figure}
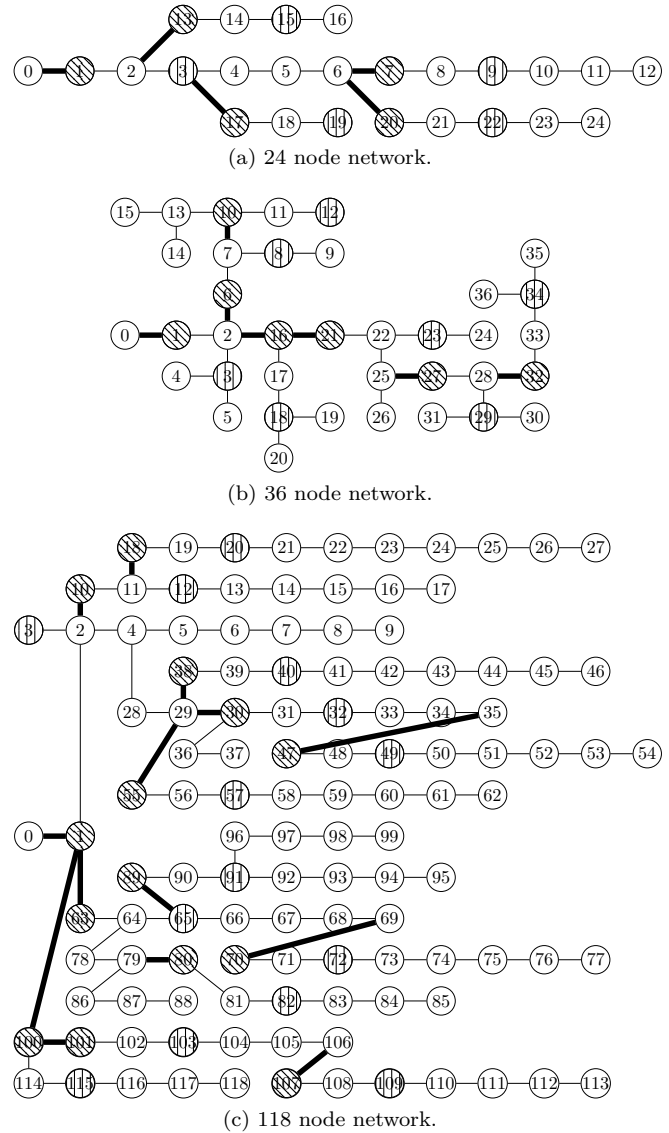

\fi

\end{document}